\documentclass{ACmod}

\DeclareFontFamily{U}{rsf}{}
\DeclareFontShape{U}{rsf}{m}{n}{
  <5> <6> rsfs5 <7> <8> <9> rsfs7 <10-> rsfs10}{}
\DeclareMathAlphabet{\mathscr}{U}{rsf}{m}{n}

\DeclareMathAlphabet{\mathgth}{U}{euf}{m}{n}

\DeclareFontFamily{U}{cyr}{}
\DeclareFontShape{U}{cyr}{m}{n}{
  <5> wncyr5 <6> wncyr6 <7> wncyr7 <8> wncyr8 <9> wncyr9 <10-> wncyr10}{}
\DeclareMathAlphabet{\mathcyr}{U}{cyr}{m}{n}

\input cyracc.def

\DeclareSymbolFont{bbold}{U}{bbold}{m}{n}
\DeclareSymbolFontAlphabet{\mathbbold}{bbold}

\makeatletter
\def\operator@font{\sf}
\makeatother

\newcommand{\Leg}{{\mathsf{Leg}}}

\newcommand{\bone}{{\mathbf 1}}

\DeclareMathOperator{\id}{id}

\DeclareMathOperator{\sym}{Sym}

\newcommand{\ra}{\rightarrow}

\renewcommand{\phi}{\varphi}

\usepackage{tikz}
\usetikzlibrary{positioning}
\usetikzlibrary{cd}
\usepackage{amscd}
\usetikzlibrary{arrows,shapes,positioning}
\usetikzlibrary{decorations.markings}
\tikzstyle arrowstyle=[scale=1]
\tikzstyle directed=[postaction={decorate,decoration={markings,
    mark=at position .65 with {\arrow[arrowstyle]{stealth}}}}]
\tikzstyle reverse directed=[postaction={decorate,decoration={markings,
    mark=at position .65 with {\arrowreversed[arrowstyle]{stealth};}}}]
\begin{document}

\title{Categorical enumerative invariants of the ground field}

\author[Junwu Tu]{% 
Junwu Tu}

\address{Junwu Tu, Institute of Mathematical Sciences, ShanghaiTech University,
Shanghai, 201210, China.} 

\begin{abstract} 
  {Abstract: For an $S^1$-framed modular operad $P$, we introduce its ``Feynman compactification" denoted by $FP$ which is a modular operad.  Let $\{\mathbb{M}^{\sf fr}(g,n)\}_{(g,n)}$ be the $S^1$-framed modular operad defined using moduli spaces of smooth curves with framings along punctures. We prove that the homology operad of $F\mathbb{M}^{\sf fr}$ is isomorphic to $H_*(\overline{M})$, the homology operad of the Deligne-Mumford operad. Using this isomorphism, we obtain an explicit formula of the fundamental class of $[\overline{M}_{g,n}/S_n]$ in terms of Sen-Zwiebach's string vertices. As an immediate application, we prove Costello's categorical enumerative invariants of the ground field match with the Gromov-Witten invariants of a point. }
\end{abstract}

\maketitle

%\tableofcontents

\section{Introduction}

\paragraph{{\bf CEI of the algebra $\mathbb{Q}$.}} Categorical enumerative invariants (CEI) were introduced by Costello~\cite{Cos2} and recently made explicit in~\cite{CCT} and~\cite{CT}. These are invariants associated with a pair $(A,s)$ where
\begin{itemize}
\item  $A$ is a cyclic $A_\infty$ algebra of dimension $d$ over a field $\mathbb{K}$ (of characteristic zero) that is proper, smooth and satisfies the Hodge-to-de-Rham degeneration.
\item $s$ is a choice of splitting of the non-commutative Hodge filtration of $A$. 
\end{itemize}
Explicitly, for each pair of non-negative integers $(g,n)$ such that $2g-2+n>0$, the CEI associated with $(A,s)$ yields a scalar
\[ \langle \alpha_1\psi_1^{k_1},\cdots,\alpha_n\psi_n^{k_n} \rangle_{g}^{A,s} \in \mathbb{K},\]
with $\alpha_1,\ldots,\alpha_n\in HH_{*-d}(A)$ in the $d$-shifted Hochschild homology of $A$, and the formal variables $\psi_1,\ldots,\psi_n$ are of homological degree $-2$. It was proved in~\cite{CT} that CEI can only be non-zero provided that 
\begin{equation}\label{eq:dim}
 \sum_{j=1}^n (-|\alpha_j|+2k_j)= 2(g-1)(3-d)+2n.
 \end{equation}
Observe that if in the cohomological degree convention, the condition above is precisely the dimension axiom in Gromov-Witten theory. Indeed, the original motivation to define and study such invariants came from mirror symmetry following Kontsevich's proposal~\cite{Kon}. It is expected that these invariants, when applied to Fukaya categories, should recover Gromov-Witten invariants at all genera. 

At first glance, comparing these categorical invariants with Gromov-Witten invariants may look rather difficult since CEI are defined purely algebraically while the GW invariants involves moduli spaces of stable maps. The purpose of this paper is to demonstrate otherwise. We shall deal with the simplest possible case when the ground field $\mathbb{K}=\mathbb{Q}$ and the algebra $A$ is also $\mathbb{Q}$. In this case, the dimension $d$ is zero, and the algebra $A$ has a unique homogeneous splitting of the Hodge filtration which we still denote by $s$. Furthermore, its Hochschild homology group is one-dimensional in degree zero $HH_*(\mathbb{Q})=\mathbb{Q}$, which implies the insertions in CEI lie in the space $\mathbb{Q}[\psi]$.  

\medskip
\noindent {{\bf Theorem A. }} {{\sl The categorical enumerative invariants of the algebra $\mathbb{Q}$ (the field of rationals)  matches the Gromov-Witten invariants of a point, i.e. we have}}
\[ \langle \psi_1^{k_1},\ldots,\psi_n^{k_n}\rangle^{\mathbb{Q},s}_{g,n} = \int_{[\overline{M}_{g,n}]} \psi_1^{k_1}\cdots\psi_n^{k_n}\]
{{\sl where on the right hand side the $\psi_j$'s represent the $\psi$-classes on the Deligne-Mumford moduli space $\overline{M}_{g,n}$. By the dimension formula~\eqref{eq:dim}, these invariants are only non-zero when $k_1+\cdots+k_n=3g-3+n$.}}

%\begin{remark}
%It is known by~\cite[Proposition 3.11]{CT} that CEI is compatible with the Givental group action. Using Givental-Teleman's reconstruction theorem~\cite{Giv,PPZ,Tel}, under mild assumptions~\cite[Assumption 5.8]{AmoTu} on the Fukaya category, the above theorem implies that CEI match with the Gromov-Witten invariants when the Fukaya category has semi-simple Hochschild cohomology. We refer to Corollary~\ref{cor:semi-simple} for a precise statement.
%\end{remark}

\paragraph{{\bf Basic constructions and notations.}}\label{para:basic} In order to prove the above theorem, it is necessary to have a better understanding of Sen-Zwiebach's string vertices~\cite{SenZwi}, since they play a central role in the construction of CEI. 

We need to introduce some notations. For a topological space $X$, denote by $C_*(X)$ its normalized singular chain complex with coefficients in $\mathbb{Q}$. For each pair $(g,n)\in \mathbb{N}\times \mathbb{N}$ such that $2g-2+n>0$, denote by $M_{g,n}^{\sf fr}$ the moduli space of tuples $(\Sigma, p_1,\ldots,p_n, \phi_1,\ldots,\phi_n)$ where $\Sigma$ is a Riemann surface of genus $g$; $p_1,\ldots,p_n$ are $n$ marked points in $\Sigma$; and for each $1\leq i\leq n$, $\phi_i: \mathbb{D}^2 \ra U(p_i)\subset \Sigma$ is a holomorphic chart such that it extends to an open neighborhood of the unit disk $\mathbb{D}^2\subset \mathbb{C}$. We also require that the pair-wise intersections of the closures of the framed disks be empty, i.e. $\overline{U(p_i)}\cap \overline{U(p_j)} =\emptyset,\;\; \forall 1\leq i<j \leq n$. As in Segal~\cite{Segal}, sewing along the coordinate charts defines two types of composition maps:
\begin{align}\label{segal-sewing}
\begin{split}
& M_{g,n}^{\sf fr} \times M_{g',n'}^{\sf fr} \ra M_{g+g',n+n'-2}^{\sf fr}, \;\;\;1\leq i\leq n, 1\leq j\leq n',\\
& M_{g,n+2}^{\sf fr} \ra M^{\sf fr}_{g+1,n}, \;\;\; 1\leq i < j \leq n+2.
\end{split}
\end{align}
We shall denote both of them by $c_{ij}$, as no confusion can arise in this way. The collection of topological spaces $\{ M_{g,n}^{\sf fr} \}_{(g,n)}$ together with the two sewing operations form a modular operad defined by Getzler-Kapranov~\cite{GetKap} in the category of topological spaces. We shall sketch a proof of this in Section~\ref{para:ex2} using Segal's argument~\cite[Section 2]{Segal}. 

Applying the normalized singular chain functor $C_*$ (with rational coefficients) to $\{ M_{g,n}^{\sf fr} \}_{(g,n)}$ yields a differential graded modular operad $\{ C_*(M_{g,n}^{\sf fr}) \}_{(g,n)}$. For the first type composition, note that the functor $C_*$ is lax monoidal, i.e. for two topological spaces $X$ and $Y$, we have a naturally defined map $C_*(X)\otimes C_*(Y) \ra C_*(X\times Y)$, the Eilenberg-Zilber map. This enables us to define a composition map still denoted by
\[ c_{ij}: C_*(M_{g,n}^{\sf fr}) \otimes C_*(M_{g',n'}^{\sf fr}) \stackrel{{\sf EZ}}{\longrightarrow} C_*(M_{g,n}^{\sf fr} \times M_{g',n'}^{\sf fr} ) \ra C_*(M_{g+g',n+n'-2}^{\sf fr}).\]
Denote this differential graded modular operad by $\mathbb{M}^{\sf fr}$, with $\mathbb{M}^{\sf fr}(g,n):= C_*(M_{g,n}^{\sf fr})$.

The modular operad $\mathbb{M}^{\sf fr}$ is in fact an {\sl $hS^1$-equivariant modular operad} with the circle actions given by rotations of the local coordinate charts. We refer to Section~\ref{sec:operads} for more details on the definition of $hS^1$-equivariant modular operad and examples. Explicitly, the circle actions on $\mathbb{M}^{\sf fr}$ are given by degree one operators
\[ B_i: \mathbb{M}^{\sf fr}(g,n) \ra \mathbb{M}^{\sf fr}(g,n), \;\; 1\leq i\leq n.\]
Denote the homotopy quotient complex of the $(hS^1)^n$-action by
\[ \mathbb{M}^{\sf fr}_{S^1}(g,n):=  \big( \mathbb{M}^{\sf fr}(g,n)[u_1^{-1},\ldots,u_n^{-1}], \partial + \sum_{i=1}^n u_i\cdot B_i\big)\]
with $u_1,\ldots,u_n$ circle parameters of homological degree $-2$. Note that the collection $\{\mathbb{M}^{\sf fr}_{S^1}(g,n)\}_{(g,n)}$ is no longer a modular operad since the the contraction map $c_{ij}$ does not descend to the quotient complexes. The symmetric group $S_n$ still acts on the equivariant complex $ \mathbb{M}^{\sf fr}_{S^1}(g,n)$. Denote its quotient complex by $\mathbb{M}^{\sf fr}_{S^1}(g,n)_{S_n}$.

\paragraph{{\bf String vertices.}} With these preparations, we are ready to introduce one of the main constructions in defining CEI: Sen-Zwiebach's differential graded Lie algebra (DGLA)~\cite[Section 3]{CCT}. As a graded vector space,  this DGLA is given by
\[     \mathfrak{g} := \big(\bigoplus_{(g,n)} \mathbb{M}^{\sf fr}_{S^1}(g,n)_{S_n} [1]\big)[[\hbar,\lambda]]\]
where  $\hbar$, $\lambda$ are two formal variables both of homological degree $-2$. The differential of $\mathfrak{g}$ is of the form $\partial + \sum_{i=1}^n u_i\cdot B_i+\hbar \Delta$ with the operator $\Delta$ defined by all possible ways of ``twisted sewing":
\[ \Delta ( \alpha \cdot u_1^{-k_1}\cdots u_n^{-k_n}) := 
 \sum_{1\leq i< j \leq n} \delta_{k_i=0}\cdot \delta_{k_j=0} \cdot c_{ij}(B_i \alpha)\]
The Lie bracket $\{-,-\}$ of $\mathfrak{g}$ is defined in a similar way:
\[ \{ \alpha \cdot u_1^{-k_1}\cdots u_n^{-k_n}, \beta \cdot v_1^{-l_1}\cdots v_{n'}^{-l_{n'}}\} := \sum_{1\leq i\leq n, 1\leq j\leq n'} \delta_{k_i=0}\cdot \delta_{l_j=0}\cdot c_{ij}(B_i \alpha, \beta)\]
Observe that the circle operator $B_i$ has homological degree one. This explains the shift by $[1]$ in the definition of $\mathfrak{g}$: to make the Lie bracket of homological degree zero. We refer to~\cite[Section 3]{CCT} for more details of the construction of $\mathfrak{g}$. 

By definition, the string vertex is a Maurer-Cartan element of the form $$\mathcal{V}:= \sum_{(g,n)} \mathcal{V}_{g,n} \hbar^g \lambda^{2g-2+n}$$ in the DGLA $\mathfrak{g}$ satisfying $\mathcal{V}_{0,3}= \frac{1}{6}{\sf pt}$ (${\sf pt}$ represents any point class in $\mathbb{M}^{\sf fr}_{S^1}(0,3)_{S_3}$. The Maurer-Cartan equation  satisfied by $\mathcal{V}$ is equivalent to the following system of equations usually known as the quantum master equation:
\begin{equation}\label{eq:qme}
(\partial+ \sum_{i=1}^n u_i\cdot B_i) \mathcal{V}_{g,n} +\Delta \mathcal{V}_{g-1,n+2} + \frac{1}{2} \sum_{g'+g''=g, n'+n''=n+2} \{ \mathcal{V}_{g',n'}, \mathcal{V}_{g'',n''}\}=0
\end{equation} 
Costello proved in~\cite[Theorem 1]{Cos2} that this system of equations together with $\mathcal{V}_{0,3}= \frac{1}{6}{\sf pt}$ has a unique solution up to gauge equivalence in $\mathfrak{g}$, which shows that the string vertex is well-defined up to gauge equivalences.

In the homological convention, a Maurer-Cartan element has degree $-1$, which implies that $\mathcal{V}_{g,n}$ is of degree $6g-6+2n$ inside the equivariant chain complex $\mathbb{M}^{\sf fr}(g,n)[u_1^{-1},\ldots,u_n^{-1}]$. Thus explicitly,  the string vertex at position $(g,n)$ is of the form
\[ \mathcal{V}_{g,n}= \sum_{(k_1,\ldots,k_n)} \mathcal{V}_{g,n}^{k_1,\ldots,k_n} u_1^{-k_1}\cdots u_n^{-k_n},\]
with the coefficient chains 
\[ \mathcal{V}_{g,n}^{k_1,\ldots,k_n} \in C_{6g-6+2n-2k_1-\cdots - 2k_n}(M_{g,n}^{\sf fr}).\]
By definition, it is also symmetric under permutations of the indices. Then, after unwinding the definition of CEI in~\cite{CCT,CT}, for the setup in the case of Theorem A, the CEI can be computed by the following formula
\[ \langle u^{k_1},\ldots,u^{k_n}\rangle^{\mathbb{Q},s}_{g,n} = \begin{cases}
0, &\;\;\; k_1+\cdots+k_n\neq 3g-3+n,\\
n! \cdot |  \mathcal{V}_{g,n}^{k_1,\ldots,k_n} |, &\;\;\; k_1+\cdots+k_n= 3g-3+n.
\end{cases}\]
Here in the second case, note that $ \mathcal{V}_{g,n}^{k_1,\ldots,k_n} $ is a zero chain, hence it is a linear combination of points. The notation $|\cdot|$ stands for the sum of its coefficients. To prove Theorem A, it remains to relate these numbers with integrals of $\psi$-classes on the Deligne-Mumford compactification $\overline{M}_{g,n}$.

\paragraph{{\bf Feynman compactification of $\mathbb{M}^{\sf fr}$.}} A key ingredient in the proof of Theorem A, introdued in Section~\ref{sec:operads}, is the notion of the {\sl Feynman compactification} of the $hS^1$-equivariant modular operad $\mathbb{M}^{\sf fr}$.  This construction yields a modular operad which we denote by $F\mathbb{M}^{\sf fr}$. Intuitively speaking, this is an ordinary modular operad such that the circle actions on the composition maps of $\mathbb{M}^{\sf fr}$ have been universally trivialized. More precisely, an element of $F\mathbb{M}^{\sf fr}(g,n)$ is given by a stable graph $G\in \Gamma((g,n))$ together with decorations: 
\begin{itemize}
\item At a vertex $v$, the decoration is by an element in the homotopy quotient $\mathbb{M}_{S^1}^{\sf fr}(g(v),{\sf Leg}(v)):= \mathbb{M}^{\sf fr}(g(v),{\sf Leg}(v))_{(S^1)^{|{\sf Leg}(v)|}}$.
\item At an edge $e$, it is decorated by a homological degree $2$ element $D_e$.
\end{itemize}
We think of the edge decoration $D_e$ as giving a universal trivialization of the circle action on the composition map, see Equation~\eqref{eq:delta} for a precise definition. 

\medskip
\noindent {{\bf Theorem B.}} {{\sl  There is an isomorphism of modular operads $$H_*(F\mathbb{M}^{\sf fr}) \cong H_*(\overline{{M}},\mathbb{Q})$$ where $\overline{{M}}$ is the Deligne-Mumford modular operad~\cite[Section 6.2]{GetKap}. }}

\begin{remark}
This result may be viewed as a higher genus extension of the beautiful works by Dotsenko-Shadrin-Vallette~\cite{DSV}~\cite{DSV2}, Drummond-Cole~\cite{Dru} and Khoroshkin-Markarian-Shadrin~\cite{KMS} in genus zero. Indeed, we prove in Section~\ref{sec:operads} that a $\mathbb{M}^{\sf fr}$-algebra $V$ together with a trivialization of the circle action induces a $F\mathbb{M}^{\sf fr}$-algebra structure on the homotopy quotient $V_{S^1}$ (see Proposition~\ref{prop:extension} for details). Together with the theorem above, this implies that the homology $H_*(V_{S^1})$ carries a $H_*(\overline{{M}})$-algebra structure. However, the approach taken here which uses the Feynman compactification $F\mathbb{M}^{\sf}$ is different even in genus zero. It is worthwhile to explore the comparison with the previous works.
\end{remark}

\begin{remark}
 We also point out that "Feynman compactifications" are {\sl different} from Getzler-Kapranov's  "Feynman transforms" of modular operads~\cite{GetKap}. Another rather interesting comparison question is to clarify the relationship between Theorem B with Getzler-Kapranov~\cite[Proposition 6.11]{GetKap}.
\end{remark}

At this point, we say a word about the proof of Theorem B. The idea is to exhibit a sequence of quasi-isomorphisms:
\[\begin{CD}
F{\mathbb{M}}^{\sf fr}(g,n) @>i^\sharp>> {\sf Tot}\big( F{\mathbb{M}}^{\sf fr}(g,n)_\bullet\big) @>\mathbb{I}>>  {\sf Tot}\big( F\widehat{\mathbb{M}}^{\sf fr}(g,n)_\bullet\big) @> p >> C_*(\overline{M}_{g,n})   
\end{CD}\]
Here $F{\mathbb{M}}^{\sf fr}(g,n)_\bullet$ is a simplicial resolution of $F{\mathbb{M}}^{\sf fr}(g,n)$ which we refer to as {{\sl Mondello's resolution}} as it's first constructed by Mondello~\cite{Mon}. The hatted version $\widehat{\mathbb{M}}^{\sf fr}(g,n)$ is given by the normalized singular chain complex of the space $\widehat{M}_{g,n}^{\sf fr}$ constructed in~\cite{KSV}. Basically, this is the framed version of the oriented real blowup of $\overline{M}_{g,n}$ along its boundary divisors, see Section~\ref{sec:operads} for more details. The first map $i^\sharp$ is constructed explicitly in Section~\ref{sec:mondello}, see Equation~\eqref{eq:in-pro}. Our main construction is the second map $\mathbb{I}$ defined in Section~\ref{sec:isomorphism}. The natural inclusion map $M_{g,n}^{\sf fr} \hookrightarrow \widehat{M}_{g,n}^{\sf fr}$, being a homotopy equivalence, does not respect the operad structures in the strict sense. The map $\mathbb{I}$ is precisely to interpolate between the two operadic compositions using higher coherent homotopies. The last map $p$ is a canonical projection map, see Paragraph~\ref{para:composition}. We prove in Section~\ref{sec:isomorphism} that the composition $p\mathbb{I}i^\sharp$ induces an isomorphism in homology that also respect the operadic composition on both sides.

\paragraph{{\bf Geometric understanding of string vertices.}} Using Theorem B, we can finally unlock the mystery of string vertices and the quantum master equation~\eqref{eq:qme}: there is a distinguished element in $H_{6g-6+2n}\big(F\mathbb{M}^{\sf fr}(g,n)\big)$ given by the homology class of following Feymann type sum
\[\sum_{G\in \Gamma((g,n))} \frac{1}{|{\sf Aut}(G)|} \prod_{e\in E_G} D_e \otimes \prod_{v\in V_G} \mathcal{V}^{\sym}_{g(v),|{\sf Leg}(v)|},\]
where the summation is over isomorphism classes of stable graphs of type $(g,n)$, $D_e$ is a degree $2$ element associated to an edge of $G$, and the superscript $\sym$ denotes the symmetrization map. Note that the closedness of this element follows from Equation~\eqref{eq:qme}. In Section~\ref{sec:fundamental}, we show that under the isomorphism $p\mathbb{I}i^\sharp$ in Theorem B, this distinguished element is mapped to the orbifold fundamental class of the symmetric Deligne-Mumford space $\overline{M}_{g,n}/S_n$, i.e. we have
\begin{equation}~\label{eq:fundamental}
(p\mathbb{I}i^\sharp)\Big( \sum_{G\in \Gamma((g,n))} \frac{1}{|{\sf Aut}(G)|} \prod_{e\in E_G} D_e \otimes \prod_{v\in V_G} \mathcal{V}^{\sym}_{g(v),|{\sf Leg}(v)|}\Big) = \frac{1}{n!}[\overline{M}_{g,n}].
\end{equation} 

Geometrically, Equation~\eqref{eq:fundamental} may be viewed as a decomposition of the fundamental class of the symmetric Deligne-Mumford space $\overline{M}_{g,n}/S_n$ according to its boundary strata labeled by stable graphs. String vertices are simply a coherent choice of the complement of a tubular neighborhood of the boundary divisor $\partial \overline{M}_{g,n}/S_n$ in the ambient space $\overline{M}_{g,n}/S_n$ for all $(g,n)$ such that $2g-2+n>0$. The coherence equation is precisely the quantum master equation in the Sen-Zwiebach DGLA. Once such a coherent choice is made, Equation~\eqref{eq:fundamental} simply corresponds to the decomposition of the fundamental class $[\overline{M}_{g,n}/S_n]$ in terms disk bundles over products of string vertices.  

We illustrate this point of view in the following figures. The left figure corresponds to the cases when $(g,n)=(0,4)$ or $(1,1)$. For example, the symmetric quotient $\overline{M}_{0,4}/S_4$ is topologically a sphere. In it there is a point (denoted by $\infty$) representing a nodal sphere with $4$ marked points. Thus the fundamental class is decomposed into a (closed) disk bundle $D_\infty$ over the point $\infty$ and $\mathcal{V}_{0,4}$ which is the closure of $\overline{M}_{0,4}/S_4 - D_\infty$. The same decomposition holds in the case of $\overline{M}_{1,1}$ with the point $\infty$ given by the once punctured nodal elliptic curve. The figure on the right illustrates the situation when two boundary divisors cross, producing a codimension $2$ strata in the moduli space. For example, the symmetric quotient $\overline{M}_{0,5}/S_5$ is decomposed into $\mathcal{V}_{0,5}$, disk bundles and double disk bundles. 
\begin{center}
\includegraphics[scale=.4]{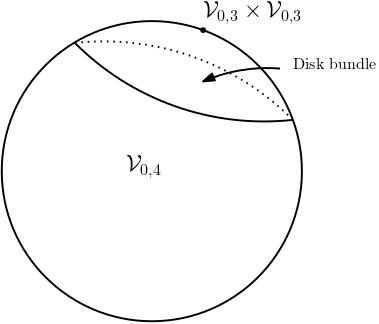}\;\;\;
\includegraphics[scale=.4]{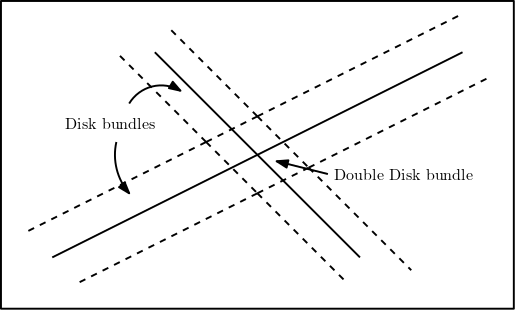}
\end{center}

This geometric point of view was already implicit in the original work of Sen-Zwiebach~\cite{SenZwi} and more recently the work of Costello-Zwiebach~\cite{CosZwi}. From this perspective, it is appropriate to think of the CEI at $(g,n)$ as a way of decomposing an integral over $\overline{M}_{g,n}$ into integrals over poly-disk bundles on each degenerate strata whose combinatorial types are given by stable graphs in $\Gamma((g,n))$.

Using the above identity of the fundamental classes, Theorem $A$ follows immediately after unwinding the definitions of CEI.

\paragraph{{\bf Organization of the paper.}} In Section~\ref{sec:operads} we define Feynman compactifications of $hS^1$-equivariant modular operads. In Section~\ref{sec:mondello} we introduce a simplicial resolution of the Feynman compactification $F\mathbb{M}^{\sf fr}$, which is a minor modification of the one constructed by Mondello~\cite{Mon}. In Section~\ref{sec:isomorphism}, we prove Theorem B. Section~\ref{sec:fundamental} is devoted to prove the main identity~\eqref{eq:fundamental}. In Section~\ref{sec:cei} we prove Theorem A.

\paragraph{{\bf Notations and Conventions.}} Throughout the paper, we work with chain complexes and homological degree convention. Unless otherwise stated, we take the ground field to be $\mathbb{Q}$, the field of rationals. For a topological space $X$, the notation $C_*(X)$ stands for its normalized singular chain complex, with coefficients in $\mathbb{Q}$.

\paragraph{{\bf Acknowledgments.}} The author is grateful for Lino Amorim, Andrei Caldararu, and Kevin Costello for discussions around the topic. The author is partially supported by the National Key Research and Development Program of China No. 2023YFA1009803 and the NSFC grant 12071290.

\section{Feynman compactifications}~\label{sec:operads}

In this section, we introduce the notion of {{\sl Feynman compactifications}} of $hS^1$-equivariant modular operads in the category of chain complexes of vector spaces over $\mathbb{Q}$. Intuitively speaking, the construction yields an ordinary modular operad such that the circle actions on the composition maps of the original operad have been universally trivialized. 

\paragraph{{\bf Recollections of modular operads.}} Let us first recall the notion of {\em modular operads} following Getzler-Kapranov~\cite[Section 2]{GetKap}.  To fix the notations, for a connected graph $G$, let us denote by
\begin{itemize}
\item $V_{G}$: the set of vertices of $G$,
\item $E_{G}$: the set of edges of $G$,
\item $\Leg(v)$: the set of half-edges at a vertex $v\in V_G$,
\item ${\sf Leaf}(G)$: the set of leaves of $G$.
\end{itemize}
A labeled connected graph is a connected graph $G$ together with a genus labeling map $g: V_G\ra \mathbb{N}$. The genus of a labeled graph is defined by $$g(G):=\sum_{v\in V_G} g(v)+\dim H^1(G).$$ 
A labeled connected graph is called a stable graph if for each $v\in V_G$ we have $2g(v)-2+n(v)>0$ where $n(v):=|\Leg(v)|$. 

Let $G$ be a stable graph. Let $e\in E_G$ be an edge of $G$. We may obtain another stable graph $G/e$ by contracting the edge $e$. Depending on $e$ being a loop edge or a non-loop edge, the contraction construction is illustrated in Figure~\eqref{figure:loop} and Figure~\eqref{figure:non-loop} respectively. 
\begin{equation}\label{figure:loop}
\begin{tikzpicture}[baseline={([yshift=-.4ex]current bounding box.center)},scale=0.5]
\draw [thick] (3.2,0) to (5.4,0);
\draw [thick] (2,0) to (3.2,0);
\draw [thick] (3.2,0) to (2,-2);
\draw [thick] (3.2,0) to (2,2);
\draw [thick] (5.2,0) circle (2); 
\node at (3.2,0) {$\bullet$};
\draw (3,0) node[label=above:{$g$}] {};
\end{tikzpicture} \;\;\;\longrightarrow\;\;\; \begin{tikzpicture}[baseline={([yshift=-.4ex]current bounding box.center)},scale=0.5]
\draw [thick] (3.2,0) to (5.4,0);
\draw [thick] (2,0) to (3.2,0);
\draw [thick] (3.2,0) to (2,-2);
\draw [thick] (3.2,0) to (2,2);
\node at (3.2,0) {$\bullet$};
\draw (4,0) node[label=above:{$g+1$}] {};
\end{tikzpicture}\end{equation}
\begin{equation}\label{figure:non-loop}
\begin{tikzpicture}[baseline={([yshift=-.4ex]current bounding box.center)},scale=0.5]
\draw [thick] (3.2,0) to (5.4,0);
\draw [thick] (2,0) to (3.2,0);
\draw [thick] (3.2,0) to (2,-2);
\draw [thick] (3.2,0) to (2,2);
\node at (3.2,0) {$\bullet$};
\node at (5.4,0) {$\bullet$};
\draw [thick] (7.2,0) to (5.4,0);
\draw [thick] (7.2,2) to (5.4,0);
\draw [thick] (7.2,-2) to (5.4,0);
\draw (3.4,0) node[label=above:{$g_1$}] {};
\draw (5,0) node[label=above:{$g_2$}] {};
\end{tikzpicture} \;\;\;\longrightarrow\;\;\; \begin{tikzpicture}[baseline={([yshift=-.4ex]current bounding box.center)},scale=0.5]
\draw [thick] (3.2,0) to (5.4,0);
\draw [thick] (2,0) to (3.2,0);
\draw [thick] (3.2,0) to (2,-2);
\draw [thick] (3.2,0) to (2,2);
\draw [thick] (3.2,0) to (5.2,-2);
\draw [thick] (3.2,0) to (5.2,2);
\node at (3.2,0) {$\bullet$};
\draw (3.5,1) node[label=above:{{\small $g_1+g_2$}}] {};
\end{tikzpicture}\end{equation}
More generally, if $I\subset E_G$ is a subset of edges of $G$, we may obtain a stable graph $G/I$ by contracting all the edges in $I$. Denote such a contraction map by
\begin{equation}
    h: G \ra G/I.
\end{equation}
In general, a morphism $h: G\ra G'$ between any two stable graphs is always given by such a contraction, up to graph isomorphisms. We refer to Getzler-Kapranov~\cite[Section 2.13-2.14]{GetKap} for a more formal treatment. Denote by $E_{G\ra G'}\subset E_G$ the set of edges that are contracted by $h$, and denote by $E_{G\ra G'}^c=E_G \backslash E_{G\ra G'}$ its complement. 

Fix a pair of non-negative integers $(g,n)$ such that $2g-2+n>0$. Following~\cite[Section 2.15]{GetKap}, let $\Gamma((g,n))$ be the category of pairs $(G,\rho)$ where $G$ is a stable graph of genus $g$ and with $n$ leaves, and $\rho: \{1,\ldots,n\}\to {\sf Leaf}(G)$ is a bijection labeling the leaves of $G$. A morphism in $\Gamma((g,n))$ is given by a morphism $h: G\ra G'$ of stable graphs which preserves the leaf labeling.

Recall from~\cite{GetKap} the notion of a stable $\mathbb{S}$-module given by a collection $\{ P(g,n)| n,g\geq 0, 2g+n-2<0\}$ of chain complexes with an action of the symmetric group $S_n$ on   $P(g,n)$. For a finite set $I$ of cardinality $n$, we set
\[ P(g,I):= \big( \bigoplus_{f:\{1,\ldots,n\} \ra I} P(g,n)\big)_{S_n},\]
with the direct sum taken over all bijections from $\{1,\ldots,n\}$ to $I$. For each stable graph $G$, denote by
\begin{equation}\label{eq:P(G)}
P(G):= \bigotimes_{v\in V_G} P(g(v), \Leg(v)).
\end{equation}
A modular operad, according to~\cite[Proposition 2.23]{GetKap}, is a stable $\mathbb{S}$-module $P$ together with an extension of $P$ (originally only defined on objects) to a functor on the category of stable graphs, that is, for each morphism $h: G\ra G'$ in $\Gamma((g,n))$ we have a morphism $P(h): P(G)\ra P(G')$ such that $P(h_1h_0)=P(h_1)P(h_0)$ whenever $h_1$ and $h_0$ are composable.

\paragraph{{\bf $hS^1$-equivariant modular operads.}} We shall work with modular operads with an $hS^1$-equivariant structure. First, we introduce some terminology. Let $(C,\partial)$ be a chain complex. By an $hS^1$-action on $C$, we mean a homological degree one operator $B: C_\bullet \ra C_{\bullet+1}$ such that $B^2=0$ and $\partial B+B\partial=0$. Let $n\geq 1$, by an $(hS^1)^n$-action, we mean $n$ {\sl mutually commuting} $hS^1$-actions $B_j\; (j=1,\ldots,n)$, such that $B_iB_j+B_jB_i=0$ for all $ 1\leq i, j\leq n$. Its homotopy quotient is defined as the chain complex
\[C_{(hS^1)^n}:= \Big( C [ u_1^{-1},\ldots, u_n^{-1} ], \partial + \sum_{j=1}^n B_j u_j\Big)\]
with the $u_j$'s of homological degree $-2$. The $\mathbb{Q}[u_1,\ldots,u_n]$-module structure is given by the quotient of $C[u_1^{\pm},\ldots,u_n^\pm]$ by the subspace $\sum_j u_j \cdot C[u_1,\ldots,u_n]$. Observe that the inclusion map
\[ C \hookrightarrow C[ u_1^{-1},\ldots, u_n^{-1} ]=C_{(hS^1)^n}\]
is a map of chain complexes, which we refer to as the canonical quotient map.

With this preparation, an $hS^1$-equivariant modular operad is given by a modular operad $P$ together with the following additional data: 
\begin{itemize}
\item There is an action of $S_n\rtimes (hS^1)^n$ on each $P(g,n)$, extending the $S_n$-action.
\item The morphism $P(h): P(G)\ra P(G')$ associated with a morphism $h: G\ra G'$ with $G,G'\in \Gamma((g,n))$ is $S_n\rtimes (hS^1)^n$-equivariant.
\item The morphism $P(h)$ is invariant under simultaneous rotation at the two ends of an edge. More precisely, we require that
\[ P(h)(B_{e,+} \alpha) = P(h) (B_{e,-} \alpha).\]
Here $\alpha \in P(G)$ is of the form $\otimes_{v\in V_G} \alpha_v$ and $B_{e,+}$ and $B_{e,-}$ are the two circle actions at the two half edges of a contracted edge $e\in E_{G\ra G'}$.
\end{itemize}

\paragraph{{\bf A first example.}}~\label{para:ex1} We give an example of $hS^1$-equivariant modular operad, essentially due to Kimura-Stasheff-Voronov~\cite{KSV}.  For each $(g,n)$, let $\widehat{M}_{g,n}$ be the oriented real blowup along the boundary divisors in the Deligne-Mumford space $\overline{M}_{g,n}$. An element in $\widehat{M}_{g,n}$ is given by a stable curve $(\Sigma, p_1,\ldots,p_n)\in \overline{M}_{g,n}$ together with a decoration at each nodal point by a unit tangent vector :
\[ v_\alpha \in T_{x_\alpha}\Sigma^n\otimes T_{y_\alpha}\Sigma^n,\]
where $\Sigma^n\ra \Sigma$ is the normalization of $\Sigma$ and $\{x_\alpha,y_\alpha\}$ is the preimage of a nodal point labeled by $\alpha$. We also consider a framed version of it denoted by $\widehat{M}_{g,n}^{\sf fr}$. An element of this moduli space $\widehat{M}_{g,n}^{\sf fr}$ consists of a decorated stable Riemann surface $(\Sigma,p_1,\ldots,p_n, \prod_\alpha v_\alpha)\in \widehat{M}_{g,n}$, together with a framing around each marked point. That is, for each $1\leq i\leq n$, a local coordinate system $\phi_i: \mathbb{D}^2 \ra U(p_i)\subset \Sigma$
such that the biholomorphic maps $\phi$'s extend to an open neighborhood of $\mathbb{D}^2$. We also require that the pair-wise intersection of the closures of the framed disks be empty, i.e. $\overline{U(p_i)}\cap \overline{U(p_j)} =\emptyset$. The collection of topological spaces $\{\widehat{M}_{g,n}\}_{(g,n)}$ forms a modular operad in the category of topological spaces.  Indeed, let $h: G\ra G'$ be a morphism of stable graphs, the composition morphism is defined as follows. The map on the underlying stable curve is the product of inclusions
\[  \prod_{v\in V_{G}}\overline{M}_{g(v),\Leg(v)}=\prod_{v'\in V_{G'}}\prod_{v\in h^{-1}(v')} \overline{M}_{g(v),\Leg(v)} \ra \prod_{v'\in V_{G'}}\overline{M}_{g(v'),\Leg(v')}\]
We keep the same framings at marked points, and the same decorations of tangent vectors at nodes corresponding to the edges of $E^c_{G\ra G'}$. For the decoration at a node corresponding to a contracted edge $e_\alpha\in E_{G\ra G'}$, we decorate by the unit tangent direction of the vector
\[ v_\alpha:= d\phi_{x_\alpha}(\frac{d}{dz})\otimes d\phi_{y_\alpha}(\frac{d}{dw})\in T_{x_\alpha}\Sigma^n\otimes T_{y_\alpha}\Sigma^n,\]
where $\phi_{x_\alpha}: \mathbb{D}_z\ra U(x_\alpha)$ and $\phi_{y_\alpha}: \mathbb{D}_w\ra U(y_\alpha)$ are the two framings at the two legs of the edge $e_\alpha$. This defines the desired composition map
\[ \xi_{G\ra G'}: \prod_{v\in V_G}\widehat{M}_{g(v),\Leg(v)}^{\sf fr}\ra\prod_{v'\in V_{G'}}\widehat{M}_{g(v'),\Leg(v')}^{\sf fr}. \]
Since the construction of $\xi_{G\ra G'}$ is local, it is straight-forward to observe the required compatibility $\xi_{G'\ra G''}\circ \xi_{G\ra G'}=\xi_{G\ra G''}$. Applying the functor $C_*(-)$ to $\{\widehat{M}_{g,n}\}_{(g,n)}$ yields a modular operad 
\[\widehat{\mathbb{M}}^{\sf fr}:= \{C_*\big(\widehat{M}_{g,n}\big)\}_{(g,n)}.\]
The $hS^1$-equivariant structure of $\widehat{\mathbb{M}}^{\sf fr}$ is defined as follows. The circle $S^1$ acts on a local coordinate system by rotation: $\phi \mapsto \phi\circ e^{i\theta}$. This action induces an action of $(S^1)^n$ on the moduli space $\widehat{M}_{g,n}$, hence an action of $(hS^1)^n$ on the chain complex $\widehat{\mathbb{M}}^{\sf fr}(g,n)=C_*\big(\widehat{M}_{g,n}\big)$. Again, it is straight-forward to verify this action is compatible with the operadic composition map defined above.

\paragraph{{\bf A second example.}}~\label{para:ex2} We give another example of $hS^1$-equivariant modular operad, essentially due to Segal~\cite{Segal}. 
For each pair $(g,n)\in \mathbb{N}\times \mathbb{N}$ such that $2g-2+n>0$, Let $M_{g,n}^{\sf fr}$ be the moduli space defined in Section~\ref{para:basic}. Analogous to the previous example, the collection of topological spaces $\{M^{\sf fr}_{g,n}\}_{(g,n)}$ forms a modular operad in the category of topological spaces. Indeed, the composition morphism
\[\eta_{G\ra G'}: \prod_{v\in V_G} M_{g(v),\Leg(v)}^{\sf fr}\ra \prod_{v\in V_{G'}} M_{g(v'),\Leg(v')}^{\sf fr}\]
associated to a contraction map $h: G \ra G'$ is defined as follows. We keep the same framing at the leaves of $G$ and $G'$. For a contracted edge $e\in E_{G\ra G'}$, we sew the Riemann surfaces at the two ends of $e$ using the two framings $\phi_{x_\alpha}: \mathbb{D}_z\ra U(x_\alpha)$ and $\phi_{y_\alpha}: \mathbb{D}_w\ra U(y_\alpha)$. More precisely, we cut out $U(x_\alpha)$ and $U(y_\alpha)$ from the two surfaces, then identify a neighborhood of the boundary circle by the equation $zw=1$. The required compatibility that $\xi_{G'\ra G''}\circ \xi_{G\ra G'}=\xi_{G\ra G''}$ is not trivial. However, as observed by Segal~\cite[Section 2]{Segal}, it readily follows from the following fact. If $\check{X}$ is a topological space obtained from a possibly disconnected Riemann surface $X$ by sewing together some of the circles making up the boundary of $X$, denote by $\pi: X\to \check{X}$ the identification map. Then the sewed surface $\check{X}$ inherits a natural complex structure defined by the property that a function $f: U \ra \mathbb{C}$ on an open subset $U\subset \check{X}$ is holomorphic if and only if the composition $f\circ \pi: \pi^{-1}(U)\to \mathbb{C}$ is holomorphic. Applying the functor $C_*(-)$ to $\{M^{\sf fr}_{g,n}\}_{(g,n)}$ yields a modular operad in the category of chain complexes:
\[\mathbb{M}^{\sf fr}:= \{C_*\big(M^{\sf fr}_{g,n}\big)\}_{(g,n)}.\]
The $hS^1$-equivariant structure of $\mathbb{M}^{\sf fr}$ is defined in the same way as in the previous example.

It is a general fact of the oriented real blowup construction that the natural inclusion map $M_{g,n}^{\sf fr} \hookrightarrow \widehat{M}_{g,n}^{\sf fr}$ is a homotopy equivalence, which induces a homotopy equivalence of chain complexes for each $(g,n)$:  
\[I_{g,n} : \mathbb{M}^{\sf fr}(g,n)\hookrightarrow \widehat{\mathbb{M}}^{\sf fr}(g,n).\]
However, observe that the operad structure of $\widehat{\mathbb{M}}^{\sf fr}$ does not restrict to that of $ \mathbb{M}^{\sf fr}$ since when composing, the former is by forming a nodal curve while the latter is by sewing using framings which yields a smooth curve.

\paragraph{{\bf Twisted sewing operations.}} Let $P$ be an $hS^1$-equivariant modular operad. We may form the homotopy quotients: for each pair $(g,n)$ we set
\[ P(g,n)_{(hS^1)^n} =  \Big( P(g,n) [ u_1^{-1},\ldots, u_n^{-1} ], \partial + \sum_{j=1}^n B_j u_j\Big)\]
This collection of stable $\mathbb{S}$-module does not form a modular operad. This is evident in the previous two examples: for instance in the case of $\mathbb{M}^{\sf fr}$, the sewing operation of two Riemann surfaces is ambiguously defined if we first quotient out the circle actions on the framings. Nevertheless, depending on sewing along a loop edge or a non-loop edge, there are two types of sewing operations by allowing a {\em full $hS^1$-twist}. We begin with the loop edge case as illustrated in Figure~\eqref{figure:loop}.

Indeed, we write down $\rho_{ij}: P(g,n)_{(hS^1)^n} \ra P(g+1, n-2)_{(hS^1)^{n-2}}$ by setting
\begin{equation}~\label{eq:self-twisted-sewing}
 \rho_{ij}(x\cdot u_1^{-k_1}\cdots u_n^{-k_n}):= \begin{cases}
0 & k_i\neq 0, \mbox{ or }\; k_j\neq 0;\\
P(h)(B_ix) u_1^{-k_1}\cdots u_n^{-k_n} & k_i=k_j=0.
\end{cases}
\end{equation}
The case of a non-loop edge contraction (Figure~\eqref{figure:non-loop}) is similar. Explicitly, we define a map of chain complexes $ \rho_{ij}: P(g,n)_{(hS^1)^n} \otimes P(h,m)_{(hS^1)^m} \ra P(g+h,n+m-2)_{(hS^1)^{n+m-2}}$ for each $1\leq i\leq n, 1\leq j\leq m$ by 
\begin{align}~\label{eq:twisted-sewing}
\begin{split}
& \rho_{ij}(x\cdot u_1^{-k_1}\cdots u_n^{-k_n}, y\cdot v_1^{-l_1}\cdots v_m^{-l_m})\\
=& \begin{cases}
0 & k_i\neq 0, \mbox{ or } l_j\neq 0;\\
P(h)(B_ix\otimes y) u_1^{-k_1}\cdots u_n^{-k_n}\cdot v_1^{-l_1}\cdots v_m^{-l_m} & k_i=l_j=0.
\end{cases}
\end{split}
\end{align}

Generalizing the above construction of twisted sewing operations, for each stable graph $G$, we put 
\[ P_{hS^1}(G) := \bigotimes_{v\in V_G} P(g(v), {\sf Leg}(v))_{(hS^1)^{|{\sf Leg}(v)|}}\]
to be the homotopy quotient version of $P(G)$, i.e. on each vertex of $G$ we put the homotopy quotient by the circle action. Furthermore, associated to an edge $e\in G$ we also have the corresponding twisted sewing map defined using $\rho_{ij}$'s constructed above:
\[ \rho_e: P_{hS^1}(G) \ra P_{hS^1}(G/e)\]
Note that $\rho_e$ is of homological degree one.

\paragraph{{ \bf Feynman compactifications of $hS^1$-equivariant modular operads.}}~\label{para:feynman} Let $P$ be an $hS^1$-equivariant modular operad. We now proceed to define its Feynman compactification $FP$ which will turn out to be a modular operad. The underlying stable $\mathbb{S}$-module of $FP$ is given by
\[ FP(g,n):= \bigoplus_{G\in {\sf G}_{g,n}} \big( D(G)\otimes P_{hS^1}(G)\big)_{{\sf Aut}(G)},\]
where ${\sf G}_{g,n}$ denotes the set of isomorphism classes of stable graphs of type $(g,n)$, and $D(G)$ is the "cocycle" (see~\cite[Section 4]{GetKap}) defined by
\[ D(G):= \bigotimes_{e\in E_G} \mathbb{Q}\cdot D_e.\]
The notation $\mathbb{Q}\cdot D_e$ stands for the one-dimensional vector space generated by  $D_e$ at homological degree $2$. Hence $D(G)$ is also one-dimensional, at homological degree $2 |E_G|$, generated by the tensor product $D_G := \bigotimes_{e\in E_G} D_e$. 
The differential on $FP(g,n)$ is the given by $\partial +\delta$ with $\partial$ the differential of $P_{hS^1}(G)$, while the extra map $\delta$ is given by
\begin{align}~\label{eq:delta}
 \delta(D_G \otimes \alpha)= \sum_{e\in E_G} D_{G/e} \otimes  \rho_e(\alpha), \;\; \alpha\in  P_{hS^1}(G).
 \end{align}
Observe that since $\rho_e$ has degree one while $D_e$ each has degree two, the degree of $\delta$ is $-1$ as required. The idea behind this definition of $\delta$ is that each $D_e$ should be thought of as a universal trivialization of the twisted sewing operation along $e$. 

\paragraph{{\bf The modular operad structure of $FP$.}}~\label{para:composition-FP} Let $G\in \Gamma((g,n))$ be a stable graph. As in Equation~\eqref{eq:P(G)}, the stable $\mathbb{S}$-module $FP$ induces a tensor product chain complex
\[ FP(G):=  \bigotimes_{v\in V_G} FP\big(g(v), \Leg(v)\big).\]
To define the modular operad structure of $FP$, we shall extend the above assignment to a functor from the category $\Gamma((g,n))$ to chain complexes. Let us consider a stable graph morphism $h: G \to G/I$ given by an edge contraction map associated with a set $I\subset E_G$. To begin with, let us consider the case when $I=\{e\}$ has only one edge. If the edge is a loop edge as in Figure~\eqref{figure:loop}, we have
\begin{align*}
 FP(G) &=  \bigoplus_{ H \in {\sf G}_{g,n+2}} \big( D(H)\otimes P_{hS^1}(H)\big)_{{\sf Aut}(H)}\\
 FP(G/e) &=  \bigoplus_{ K \in {\sf G}_{g+1,n}} \big( D(K)\otimes P_{hS^1}(K)\big)_{{\sf Aut}(K)}
\end{align*}
For a stable graph $H\in {\sf G}_{g,n+2}$, we denote by $H_{ij}\in {\sf G}_{g+1,n}$ the stable graph obtained from $H$ by forming a loop $e$ using the two legs labeled by $i$ and $j$. Observe that $D(H_{ij}/e)\cong D(H)$ since both graph have isomorphic set of edges. Then we set the composition map as
\begin{align}~\label{eq:composition1}
\begin{split}
& FP(h)\big( D_H \otimes x u_1^{-k_1}\cdots u_{n+2}^{-k_{n+2}} \big)\\
= & \begin{cases}
D_{H_{ij}} \otimes (x \cdot u_1^{-k_1}\cdots u_{n+2}^{-k_{n+2}})\cdot (u_i+u_j)    & k_i\neq 0, \mbox{ or }\; k_j\neq 0\\
D_{H_{ij}/e}\otimes P(h)(x) \cdot u_1^{-k_1}\cdots u_{n+2}^{-k_{n+2}} & k_i=k_j=0 
\end{cases}
\end{split}
\end{align}
\begin{Lemma}
The composition map $FP(h)$ is a chain map.
\end{Lemma}
\begin{proof}
For simplicity, we shall suppress the $u$'s to only $u_i$ and $u_j$ as other circle parameters are not relevant in the calculation. We need to verify that $[\partial+uB+\delta, FP(h)]=0$. When the input is of the form $D_H\otimes x$, one has
\begin{align*}
(\partial+uB+\delta) FP(h)(D_H\otimes x) & = (\partial+uB+\delta) (D_{H_{ij}/e}\otimes P(h)(x)) \\
&= \delta D_{H_{ij}/e} \otimes P(h)(x) + D_{H_{ij}/e}\otimes \partial P(h)( x),\\
FP(h) (\partial+uB+\delta)(D_H\otimes x) &= FP(h) (\delta D_H \otimes x + D_H \otimes \partial x)\\
&= \delta D_{H_{ij}/e} \otimes P(h)(x)+ D_{H_{ij}/e}\otimes P(h)(\partial x).
\end{align*}
The modular operad structure $P(h)$ is a chain map, which shows that the two computations indeed give the same answer. When the input is of the form $D_H\otimes x u_i^{-1}$ we compute the two compositions as
\begin{align*}
& (\partial+uB+\delta) FP(h)(D_H\otimes xu_i^{-1}) \\
=  & (\partial+uB+\delta) (D_{H_{ij}}\otimes x) \\
= & \delta D_{H_{ij}} \otimes x + D_{H_{ij}}\otimes \partial x\\
& FP(h) (\partial+uB+\delta)(D_H\otimes x u_i^{-1}) \\
= &  FP(h) (\delta D_H \otimes x u_i^{-1} + D_H \otimes \partial x u_i^{-1} + D_H\otimes B_i x)\\
= & \big(\delta D_{H_{ij}} \otimes x - D_{H_{ij}/e}\otimes P(h)(B_ix)\big)+ D_{H_{ij}}\otimes \partial x+D_{H_{ij}/e}\otimes P(h)(B_ix)\\
= & \delta D_{H_{ij}} \otimes x + D_{H_{ij}}\otimes \partial x
\end{align*}
The case of $D_H\otimes x u_j^{-1}$ is similar. In the general case when the input is $D_H\otimes x u_i^{-k_i}u_j^{-k_j}$ with $k_i, k_j\geq 1$, the computation is easier and is left as an exercise.
\end{proof}
The case when $e$ is a non-loop edge as in Figure~\eqref{figure:non-loop} is defined similarly:
\begin{align*}
& FP(h)\big( D_{H_1} \otimes x \cdot u_1^{-k_1}\cdots u_{n}^{-k_{n}}, D_{H_2} \otimes y\cdot v_1^{-l_1}\cdots v_m^{-l_m} \big)\\
= & \begin{cases}
D_{H_{ij}} \otimes \big( x \cdot u_1^{-k_1}\cdots u_{n}^{-k_{n}})\otimes y\cdot v_1^{-l_1}\cdots v_m^{-l_m} \big)(u_i+v_j)   & k_i\neq 0, \mbox{ or }\; l_j\neq 0\\
D_{H_{ij}/e}\otimes P(h)(x,y) \cdot u_1^{-k_1}\cdots u_{n}^{-k_{n}}\cdot v_1^{-l_1}\cdots v_m^{-l_m} & k_i=k_j=0 
\end{cases}
\end{align*}
The verification that this defines a chain map is analogous to the previous lemma, and will be omitted.

For a general $I\subset E_G$, we may factor the map $h: G\to G/I$ as
\[ G \stackrel{h_1}{\to} G/\{e_1\} \stackrel{h_2}{\to} G/\{e_1,e_2\} \to \cdots \stackrel{h_{|I|}}{\to} G/I,\]
and define $FP(h):= FP(h_1)\cdots FP(h_{|I|})$. The well-definedness of $FP(h)$, that it is independent of the factorization of $h$, is included in the proof of the following proposition. 

\begin{Proposition}
    The stable $\mathbb{S}$-module $FP$ forms a modular operad, i.e. for each stable pair $(g,n)$, the construction above defines a functor $FP: \Gamma((g,n)) \to {\sf Ch}_\mathbb{Q}$. 
\end{Proposition}

\begin{proof}
We first argue the map $FP(h)$ is well-defined. Indeed, it is a map the following form:
\[ \bigotimes_{v\in V_G} FP\big(g(v), \Leg(v)\big)= \bigotimes_{w\in V_{G/I}} \Big( \bigotimes_{v\in h^{-1}(w)} FP\big(g(v),n(v)\big)\Big) \stackrel{FP(h)}{\longrightarrow} \bigotimes_{w\in V_{G/I}} FP\big( g(w),n(w)\big).\]
For each vertex $w\in V_{G/I}$, the pre-image of a small neighborhood of $w$ under $h$ is a stable subgraph $G_w$ in $G$ of type $(g(w),n(w))$. Given an element
\[ \alpha_w \in  \bigotimes_{v\in h^{-1}(w)} FP\big(g(v),n(v)\big),\]
then the vertex decoration at $w\in G/I$ is obtained from $\alpha_w$ by applying one of the following two types of operations at each $e\in E_{G_w}\bigcap I$:
\begin{enumerate}
    \item if in $\alpha_w$ the circle parameters at the two ends of $e$ both  have powers $k_i=k_j=0$, the we contract this edge using $P$;
    \item otherwise, we keep the edge and multiple $\alpha_w$ by $u_i+u_j$.
\end{enumerate}
The well-definedness of $FP$ follows from observing that these two types operations at edges all commutate with each other. The same argument also applies to prove that $FP$ is a functor on $\Gamma((g,n))$.
\end{proof}

\paragraph{{\bf Trivializations of circle actions.}} By a cyclic $P$-algebra structure on a chain complex $(V,b)$ we mean the following data:
\begin{itemize}
\item An $hS^1$-action $B$ on the chain complex $(V,b)$.
\item An $(hS^1)^n\ltimes S_n$-equivariant action map
\[ \rho_V: P(g,n) \ra V^{\otimes n}.\]
\item A graded symmetric inner product $\langle-,-\rangle: V^{\otimes 2} \ra \mathbb{K}$ such that $b$ is graded anti-self-adjoint and $B$ is graded self-adjoint. 
\end{itemize}
The compatibility condition between $\rho_V$ and the inner product should be self-evident:  for an elementary contraction $c_{ij}$, it corresponds to contracting two copies of $V$ at the half-edges labeled by $i$ and $j$ using the inner product.

A trivialization of the circle action $B$ (see~\cite{DSV,KMS}) is given by a chain map
\[ s: (V,b) \ra (V[[u]],b+uB)\]
of the form $s(x)=x+R_1(x)u+R_2(x)u^2+\cdots$ with $R_j: V \ra V, \; (j\geq 1)$ of homological degree $2j$. 

\begin{Proposition}~\label{prop:extension}
Let $V$ be a cyclic $P$-algebra and let $s: V\ra V[[u]]$ be a trivialization of the circle action. Then this data induces a natural cyclic $FP$-algebra structure on $V_{hS^1}= V[u^{-1}]$. 
\end{Proposition}

\begin{proof}
Taking $(hS^1)^n$-quotients of the action map
\[  \rho_V: P(g,n) \ra V^{\otimes n}\]
yields an action of $P_{hS^1}(g,n)$ on $V_{S^1}$, i.e. an action map still denoted by
\[ \rho: P_{hS^1}(g,n) \ra V_{hS^1}^{\otimes n}.\] 
To extend this action to $FP$ it suffices to define the action of $D_e$ associated with an edge $e\in E_G$ of a stable graph $G$.  Observe that the inverse operator of $R$ is another operator of the form
\[ T= \id+T_1u+T_2u^2+\cdots.\]
Here the operator $T_j$'s are defined by the following identity
\[ \sum_{i+j=k} T_iS_j = \begin{cases} \id \mbox{\;\; if \;\;} k=0,\\ 0 \mbox{\;\; if \;\;} k\geq 1.\end{cases}\]
Solving the above recursively yields formulas of $T_j$ in terms $R_i$'s. For example, we have
 \[ T_1=-R_1, \;\; T_2= -R_2+R_1R_1 \]
Using the $R$'s and $T$'s, we define a contraction map
$$H:=\sum_{i\geq 0, j\geq 0} H_{i,j} : V[u^{-1}] \otimes V[u^{-1}] \ra \mathbb{K}$$   with component maps $H_{i,j}: V \cdot u^{-i} \otimes V\cdot u^{-j}  \ra \mathbb{K}$ given by formula
\[H_{i,j}( xu^{-i} , yu^{-j}  ):= \langle (-1)^{j} \sum_{l=0}^j R_{l} T_{i+j+1-l} x, y\rangle.\]
Here we set $R_0=T_0=\id$. Then we define the action of $D_e$ on $V_{hS^1}$ by the symmetrization of $H$, i.e. $D_e$ acts by
\[ H^{\sf sym} (\alpha,\beta) := \frac{1}{2} \big(H(\alpha,\beta)+(-1)^{|\alpha||\beta|}H(\beta,\alpha)\big).\]
The verification that $[b+uB,H]$ is compatible with the boundary map $\delta D_e$ in Equation~\eqref{eq:delta} is a calculation, and is done in~\cite[Proposition 4.5]{CT}.
\end{proof}

\section{Mondello's resolution}~\label{sec:mondello}

Let $\mathbb{M}^{\sf fr}$ and $\widehat{\mathbb{M}}^{\sf fr}$ be the two $hS^1$-equivariant modular operads as defined in Paragraphs~\ref{para:ex1} and~\ref{para:ex2}. In this section, we introduce a simplicial resolution denoted by $ F\widehat{\mathbb{M}}^{\sf fr}(g,n)_\bullet$ of the Feynman compactification $F\widehat{\mathbb{M}}^{\sf fr}(g,n)$. We call it Mondello's resolution since its construction is a minor modification of Mondello's construction in~\cite{Mon}. Similarly, we also have the resolution $F\mathbb{M}^{\sf fr}(g,n)_\bullet$   of the Feynman compactification $F\mathbb{M}^{\sf fr}(g,n)$. The main result of the section is the existence of the following quasi-isomorphisms:
\begin{align*}
F\mathbb{M}^{\sf fr}(g,n) &\cong {\sf Tot}\big(  F\mathbb{M}^{\sf fr}(g,n)_\bullet\big)\\
F\widehat{\mathbb{M}}^{\sf fr}(g,n) &\cong {\sf Tot}\big( F\widehat{\mathbb{M}}^{\sf fr}(g,n)_\bullet\big) 
\end{align*}
We shall mainly deal with the hatted version, the un-hatted version is completely parallel.

\paragraph{{\bf Mondello's construction~\cite{Mon}.}} Let $X$ be a smooth algebraic variety  with a simple~\footnote{Simplicity is not necessary for this construction, as was remarked by Mondello~\cite[Footnote 1]{Mon}.} normal crossing divisor
\[ D= \bigcup_{i\in I} \overline{D}_i,\]
i.e. we assume that each $\overline{D}_i$ is a smooth irreducible divisor, and pair-wise intersect transversally. For a subset $J\subset I$, denote by $\overline{D}_J:=\bigcap_{j\in J} \overline{D}_j$. Denote by $\widehat{D}_J$ the manifold with corners obtained by performing a real oriented blow-up of $\overline{D}_J$ along a divisor $\bigcup_{i\in I \backslash J} \overline{D}_{J\cup\{i\}}$. Note that the natural inclusion $D_J:=\overline{D}_J\backslash ( \bigcup_{i\in I \backslash J} \overline{D}_{J\cup\{i\}}) \hookrightarrow \widehat{D}_J$ is a homotopy equivalence. For a chain $J_0\subsetneq J_1\subsetneq \cdots \subsetneq J_k\subset I$ of inclusions of subsets of $I$, define $\widehat{D}_{J_0,\cdots,J_k} := \widehat{D}_{J_0} \times_{\overline{D}_{J_0}} \overline{D}_{J_k}$. Observe that the space 
\begin{equation}~\label{eq:circle-bundle}
 \widehat{D}_{J_0,\cdots,J_k}=\prod_{i\in J_k\backslash J_0} b_{J_k}^*\big( SN_{\overline{D}_i/X}|_{\overline{D}_{J_k}}\big)\end{equation}
is given by the $(S^1)^d$-bundle where $SN_{\overline{D}_i/X}$ is the circle bundle associated to the normal bundle of $\overline{D}_i$ in $X$, the map $b_{J_k}: \widehat{D}_{J_k} \ra \overline{D}_{J_k}$ is the natural projection, and $d=|J_k\backslash J_0|$. In {\sl loc. cit.} the author constructs a simplicial topological space $X_\bullet$ with
\[ X_k:= \coprod_{J_0\subsetneq J_1\subsetneq \cdots \subsetneq J_k\subset I} \widehat{D}_{J_0,\ldots,J_k},\]
such that its geometric realization $|X_\bullet|$ is homotopy equivalent to $X$. 

\paragraph{{\bf Application to $\overline{M}_{g,n}$.}} We apply the previous construction to the case with $X$ equals the Deligne-Mumford compactification $\overline{M}_{g,n}$, and $D$ the divisor corresponding to nodal curves. Instead of working with an index set $I$ labeling the irreducible components of $D$, it is most natural in this setup to work with stable graphs in $\Gamma((g,n))$. That is, for each $G\in \Gamma((g,n))$, let us write $\overline{M}_G$ for the image of the map
\begin{equation}~\label{eq:xi}
  \xi_G: \prod_{v\in V_G^{\sf int}}\overline{M}_{g(v),|{\sf Leg}(v)|} \ra \overline{M}_G\subset \overline{M}_{g,n}.
  \end{equation}
Denote by $M_G:=\xi_G\big( \prod_{v\in V_G^{\sf int}}M_{g(v),{\sf val}(v)}\big)$ the  locus of $\overline{M_G}$ such that at each vertex of $G$ we have a smooth curve in $M_{g(v),|{\sf Leg}(v)|}$. The oriented real blowup of $\overline{M}_G$ along boundary divisors is denoted by $\widehat{M}_G$. A general fact of oriented real blowup is that the natural inclusion map $M_G\hookrightarrow \widehat{M}_G$ is a homotopy equivalence. Let $G \ra G'$ be a graph contraction map. It induces a corresponding inclusion map
\[ \xi_{G\ra G'} : \overline{M}_G \ra \overline{M}_{G'}\]

In dimension $k\geq 0$, we describe the set of non-degenerate $k$-simplices of Mondello's construction applied to $\overline{M}_{g,n}$ in terms of stable graphs. For each $k$-step graph contraction, we set
\[ \widehat{M}_{G_0\ra\cdots\ra G_k}:=\widehat{M}_{G_k}\mid_{\overline{M}_{G_0}}.\]
Explicitly, an element of $\widehat{M}_{G_k}\mid_{\overline{M}_{G_0}}$ consists of
\begin{itemize} 
\item A stable curve with marked points $(\Sigma,p_1,\ldots,p_n)\in \overline{M}_{G_0}$. By definition, its dual graph $G$ admits a contraction map $G\ra G_0$.
\item For each edge that got contracted by the composition map $G\ra G_0\ra G_k$,  the corresponding node $x\in \Sigma$ is decorated by a unit tangent vector in the tensor product $T_{x_+}\Sigma^n\otimes T_{x_-}\Sigma^n$ where $\Sigma^n\ra \Sigma$ is the normalization of the nodal curve $\Sigma$, and $\{x_+,x_-\}$ is the pre-image of the nodal point $x$.
\end{itemize}  
Then the set of non-degenerate $k$-simplices is given by
\[ \mathbb{X}_k:=\varinjlim_{G_0 \ra\cdots\ra G_k} \widehat{M}_{G_0\ra\cdots\ra G_k}\]
where the colimit is taken over the category of isomorphisms of $k$-step stable graph contractions.

Next, we describe the boundary maps of $\mathbb{X}_\bullet$. There are three types of boundary maps in $\mathbb{X}_k \ra \mathbb{X}_{k-1}$: 
\begin{itemize}
\item The first type boundary map is a map of the form
\begin{equation}\label{eq:boundary1}
 \delta_1: \widehat{M}_{G_0\ra\cdots\ra G_{k}}=\widehat{M}_{G_k}\mid_{\overline{M}_{G_0}} \; \ra \widehat{M}_{G_0\ra\cdots\ra G_{k-1}}=\widehat{M}_{G_{k-1}}\mid_{\overline{M}_{G_0}}
 \end{equation}
given by the natural projection map which forgets the decorations at nodes that got contracted in the map $G_{k-1}\ra G_k$.
\item The second type boundary map is given by
\begin{equation}\label{eq:boundary2}
 \delta_2: \widehat{M}_{G_0\ra\cdots\ra G_{k}}=\widehat{M}_{G_k}\mid_{\overline{M}_{G_0}} \; \ra \widehat{M}_{G_1\ra\cdots\ra G_{k}}=\widehat{M}_{G_k}\mid_{\overline{M}_{G_1}},
 \end{equation}
simply the inclusion map induced by the inclusion $\overline{M}_{G_0}\hookrightarrow \overline{M}_{G_1}$. 
\item The third type boundary map is defined in the case $k\geq 2$ and for each $1\leq j\leq k-1$, and is given by the identity map
\begin{equation}\label{eq:boundary3}
 \delta_3: \widehat{M}_{G_0\ra\cdots\widehat{G_j}\cdots\ra G_{k}}=\widehat{M}_{G_k}\mid_{\overline{M}_{G_0}} \; \ra \widehat{M}_{G_0\ra\cdots\widehat{G_j}\cdots\ra G_{k}}=\widehat{M}_{G_k}\mid_{\overline{M}_{G_0}},
 \end{equation}
 \end{itemize} 

\paragraph{{\bf Circle/Disk bundles.}} For each stable graph contraction $G_0\ra\cdots\ra G_k$, we shall work with a modified but homotopy equivalent version of $\widehat{M}_{G_0\ra\cdots\ra G_k}$. Let $G\ra G'$ be a graph contraction. Denote by $N_{G\ra G'}$ the normal bundle associated with the embedding $\overline{M}_G \hookrightarrow \overline{M}_{G'}$. Its rank is equal to $d=|E_{G\ra G'}|$, the number of contracted edges in $G\ra G'$. Moreover, we denote by $SN_{G\ra G'}$ the associated $(S^1)^d$-bundle, and $DN_{G\ra G'}$ the associated $D^d$-bundle. For $G \ra \star$ with $\star$ the unique stable graph in $\Gamma((g,n))$ with one vertex, we write $SN_G$ and $DN_G$ for the corresponding circle/disk bundles. As observed by Mondello~\cite{Mon} (see Equation~\ref{eq:circle-bundle}), we have
\[ \widehat{M}_{G_0\ra\cdots\ra G_k}= \widehat{M}_{G_k}\mid_{\overline{M}_{G_0}}= b_{G_0}^*(SN_{G_0\ra G_k}),\]
where $b_{G_0}: \widehat{M}_{G_0} \ra \overline{M}_{G_0}$ is the canonical projection map by definition of the oriented real blowup construction. 

We shall make two modifications to the above space, both of which are homotopy equivalent to it by construction:
\begin{itemize}
\item The first modification is to add a disk bundle in the normal direction of $\overline{M}_{G_k}$ inside $\overline{M}_{g,n}$. More precisely, we use 
\[ b_{G_0}^*\big( SN_{G_0\ra G_k} \times_{\overline{M}_{G_0}} DN_{G_k}|_{\overline{M}_{G_0}}\big).\]
Since disks are contractible, it is clear that this is homotopy equivalent to $\widehat{M}_{G_0\ra\cdots\ra G_k}$.
\item The second modification is to put framings on the marked points, then realizing the above space as its $(hS^1)^n$-quotient. Since the circle action on framings is a free action, the resulting space is also homotopy equivalent to $\widehat{M}_{G_0\ra\cdots\ra G_k}$.
\end{itemize} 

More formally, for each $k\geq 0$, we define 
$$\mathbb{Y}_k:=\varinjlim_{G_0\ra\cdots\ra G_k} \mathbb{Y}^{\sf fr}_{G_0\ra\cdots\ra G_k}/(hS^1)^n$$
where an element of $\mathbb{Y}^{\sf fr}_{G_0\ra\cdots\ra G_k}$ consists of a stable curve with marked points together with framings $(\Sigma,p_1,\ldots,p_n,\phi_1,\ldots,\phi_n)\in \overline{M}^{\sf fr}_{G_0}$. This stable curve's dual graph $G$ admits a contraction map $G\ra G_0$. For each edge that got contracted by the composition map $G\ra G_0\ra G_k$, we decorate the corresponding node by a unit tangent vector in the tensor product $T_{x_+}\Sigma^n\otimes T_{x_-}\Sigma^n$. For each of the remaining node $y$ of $\Sigma$, we decorate it by a tangent vector in the unit disk of $T_{y_+}\Sigma^n\otimes T_{y_-}\Sigma^n$.

Analogous to $\mathbb{X}_\bullet$, there are three types of simplicial boundary maps on $\mathbb{Y}_\bullet$ so that it is compatible with the boundary map of $\mathbb{X}_\bullet$ defined in Equations~\eqref{eq:boundary1}~\eqref{eq:boundary2}~\eqref{eq:boundary3}. We shall  abuse the notations $\delta_i\, (i=1,2,3)$ for the three types of boundary maps.
 \begin{itemize}
\item For the first type boundary map, we set
\[ \delta_1: \mathbb{Y}_{G_0\ra\cdots\ra G_k}^{\sf fr}/(hS^1)^n \ra \mathbb{Y}_{G_0\ra\cdots\ra G_{k-1}}^{\sf fr}/(hS^1)^n\]
to be the natural map induced by the inclusion map $S^1\subset \mathbb{D}$ at the nodes that are contracted by the morphism $G_{k-1}\ra G_{k}$.
\item The boundary map of the second type
\[\delta_2 : \mathbb{Y}_{G_0\ra\cdots\ra G_k}^{\sf fr}/(hS^1)^n \ra \mathbb{Y}_{G_1\ra\cdots\ra G_k}^{\sf fr}/(hS^1)^n\]
is again induced by the inclusion $\widehat{M}_{G_0}\ra \widehat{M}_{G_1}$ while keeping decorations at nodes and framings at marked points.
\item The third type boundary map is again the identity map:
\[ \delta_3=\id: \mathbb{Y}_{G_0\ra\cdots\ra G_k}^{\sf fr}/(hS^1)^n\ra \mathbb{Y}_{G_0\ra\cdots\widehat{G_j}\cdots\ra G_k}^{\sf fr}/(hS^1)^n.\]
\end{itemize}

Note that by construction we have a homotopy equivalence of simplicial spaces between $\mathbb{X}_\bullet$ and $\mathbb{Y}_\bullet$. By Mondello's theorem, we obtain a quasi-isomophism of chain complexes
\[  C_*(\overline{M}_{g,n}) \cong {\sf Tot} \big( C_*(\mathbb{Y}_\bullet) \big)\]

\paragraph{{\bf Simplicial resolution of the Feynman compactification.}} To link the simplicial complex $C_*(\mathbb{Y}_\bullet)$ with the Feynman compactification $F\widehat{\mathbb{M}}^{\sf fr}(g,n)$ we observe that there is a discrepancy at nodes. Indeed, let us consider a node $x$ corresponding to a contracted edge in the dual graph of a surface $\Sigma$. In the case of $\mathbb{Y}_\bullet$, a node $x$ is decorated by elements in the tensor product $T^*_{x_+}\Sigma^n\otimes T^*_{x_-}\Sigma^n$, while for $F\widehat{\mathbb{M}}^{\sf fr}(g,n)$ the two circle actions on $T^*_{x_+}$ and $T^*_{x_-}$ has already been quotiented out. The two constructions are related by the following
\begin{Lemma}~\label{lem:circle}
Let $M$ and $N$ be two complexes with a circle action. Then there is a canonical homotopy equivalence
\[ (M\otimes N)_{S^1} \cong M_{S^1}\otimes N_{S^1}\otimes \mathbb{K}[\epsilon],\]
where on the left hand side we take the homotopy quotient of the off-diagonal action, i.e. with circle operator $B_M-B_N$,  and on the right hand side, the differential is given by $b_M+u_1B_M+b_N+u_2B_N+(u_1+u_2)\epsilon$.
\end{Lemma}

\begin{proof}
This is standard Koszul duality between $\mathbb{K}[\epsilon]$ and its Koszul dual coalgebra $u^{-1}\mathbb{K}[u^{-1}]$. Thus for any $hS^1$-module $Q$, there is a chain equivalence
\[ \big( (Q)_{S^1}\otimes \mathbb{K}[\epsilon], b_Q+uB_Q+u\epsilon\big)\cong Q.\]
Apply this formula to our setting to obtain
\begin{align*}
M_{S^1}\otimes N_{S^1}\otimes \mathbb{K}[\epsilon] & \cong (M\otimes N)_{S^1\times S^1} \otimes \mathbb{K}[\epsilon]\\
&\cong \big((M\otimes N)_{S^1}\big)_{S^1}\otimes \mathbb{K}[\epsilon]\\
&\cong (M\otimes N)_{S^1}.
\end{align*}
where in the second equivalence we split the product $hS^1$-action using $B_M-B_N$ and $B_M+B_N$.
\end{proof}

We proceed to define a simplicial resolution of $F\widehat{\mathbb{M}}^{\sf fr}(g,n)$.  Let $G_0\ra\cdots\ra G_k$ be a $k$-step contraction. Define
\[ F\widehat{\mathbb{M}}^{\sf fr}_{G_0\ra\cdots\ra G_k}:= \bigotimes_{v\in V_{G_0}}\widehat{\mathbb{M}}_{S^1}^{\sf fr}(g(v),{\sf Leg}(v)) \otimes \bigotimes_{e\in E_{G_0\ra G_k}} C_*(S^1) \otimes \bigotimes_{e'\in E^c_{G_0\ra G_k}} C_*(D).\]
That is, 
\begin{itemize}
\item For each vertex we decorate it by homotopy $(S^1)^{|{\sf Leg}(v)|}$-quotient of $\widehat{\mathbb{M}}^{\sf fr}(g(v),{\sf Leg}(v))$.
\item For each contracted edge $e\in E_{G_0\ra G_k}$, we decorate by the cellular chain complex of circle:
\[C_*(S^1):=\mathbb{K}[\epsilon]=\mathbb{K}\cdot \bone \oplus \mathbb{K}\cdot \epsilon\]
with $\epsilon$ of homological degree $1$. 
\item For each of the remaining edges in $E_{G_0\ra G_k}^c$, we decorate by the cellular chain complex of the disk 
\[C_*(D):= \mathbb{K}\cdot \bone \oplus \mathbb{K}\cdot \epsilon \oplus \mathbb{K}\cdot D,\]
with $D$ of homological degree $2$. The differential acts on the $2$-cell by
$\partial(D)=\epsilon$.
\end{itemize}
The differential of $F\widehat{\mathbb{M}}^{\sf fr}_{G_0\ra\cdots\ra G_k}$  is not simply the tensor product differential. The extra differential comes from Lemma~\ref{lem:circle}.  Indeed, for every edge $e\in E_{G_0}$ (no matter if it's contracted by $G_0\ra G_k$ or not), the extra differential is given by
\[ \partial (\alpha_+\otimes \bone_e \otimes \alpha_-)=(u_{e_+} \alpha_+)\otimes \epsilon_e\otimes \alpha_-+\alpha_+\otimes \epsilon_e\otimes (u_{e_-}\alpha_-),\]
with $u_{e_+}$ and $u_{e_-}$ the two circle parameters acting on chains $\alpha_+$ and $\alpha_-$ on the two ends of $e$. For each $k\geq 0$, we define the set of $k$-simplices by setting
 \[ F\widehat{\mathbb{M}}^{\sf fr}(g,n)_k:= \varinjlim_{G_0\ra\cdots\ra G_k} F\widehat{\mathbb{M}}^{\sf fr}_{G_0\ra\cdots\ra G_k}\] 
to obtain a simplicial chain complex $F\widehat{\mathbb{M}}^{\sf fr}(g,n)_\bullet$. Again, just like $\mathbb{X}_\bullet$ and $\mathbb{Y}_\bullet$, it has  three types of simplicial boundary maps: 
\begin{itemize}
\item The first type boundary map
\[ \delta_1: (F\widehat{\mathbb{M}}^{\sf fr})_{G_0\ra\cdots\ra G_k} \ra (F\widehat{\mathbb{M}}^{\sf fr})_{G_0\ra\cdots\ra G_{k-1}}\]
is defined via the inclusion $C_*(S^1)\hookrightarrow C_*(D)$ at the edges that are contracted by $G_{k-1}\ra G_{k}$.
\item  The second type boundary map 
\[\delta_2: (F\widehat{\mathbb{M}}^{\sf fr})_{G_0\ra\cdots\ra G_k} \ra (F\widehat{\mathbb{M}}^{\sf fr})_{G_1\ra\cdots\ra G_{k}}\]
is defined by applying Lemma~\ref{lem:circle} at the edges that are contracted by the map $G_0\ra G_1$.
\item The third simplicial boundary map $$\delta_3: (F\widehat{\mathbb{M}}^{\sf fr})_{G_0\ra\cdots\ra G_k} \ra (F\widehat{\mathbb{M}}^{\sf fr})_{G_0\ra\cdots\widehat{G_j}\cdots\ra G_k}$$ is again the identity map.
\end{itemize}

In conclusion, for each pair $(g,n)$ such that $2g-2+n>0$, we have defined a simplicial chain complex $F\widehat{\mathbb{M}}^{\sf fr}(g,n)_\bullet$. Denote its total complex by
\[ {\sf Tot}\big(F\widehat{\mathbb{M}}^{\sf fr}(g,n)_\bullet\big):=\big( \bigoplus_{k\geq 0} \varinjlim_{G_0\ra\cdots\ra G_k} (F\widehat{\mathbb{M}}^{\sf fr})_{G_0\ra\cdots\ra G_k}[k], \partial+\delta\big). \]
Mondello~\cite{Mon}'s simplicial resolution construction implies the following

\begin{Theorem}~\label{thm:mondello}
There is a homotopy equivalence of chain complexes
\[ C_*(\overline{M}_{g,n})\cong {\sf Tot}\big(F\widehat{\mathbb{M}}^{\sf fr}(g,n)_\bullet\big).\]
\end{Theorem}

\paragraph{{\bf Homotopy colimit of $F\widehat{\mathbb{M}}^{\sf fr}(g,n)_\bullet$.}} When $k=0$, we have
\[ F\widehat{\mathbb{M}}^{\sf fr}(g,n)_0 = \varinjlim_{G\in \Gamma((g,n))}  \bigotimes_{v\in V_{G}}\widehat{\mathbb{M}}_{S^1}^{\sf fr}(g(v),{\sf Leg}(v)) \otimes  \bigotimes_{e\in E_{G}} C_*(D).\]
On the other hand the Feynman compactification is given by
\[ F\widehat{\mathbb{M}}^{\sf fr}(g,n)=\varinjlim_{G\in \Gamma((g,n))}  \bigotimes_{v\in V_{G}}\widehat{\mathbb{M}}_{S^1}^{\sf fr}(g(v),{\sf Leg}(v)) \otimes  \bigotimes_{e\in E_{G}} D_e \cdot \mathbb{Q}.\]
The latter space is clearly a subspace of the former. However, observe that neither the canonical inclusion map nor the canonical projection map is a map of complexes. 

In the following proposition, we prove that the homotopy colimit of the simplicial chain complex $F\widehat{\mathbb{M}}^{\sf fr}(g,n)_\bullet$ is given by the Feynman compactification $F\widehat{\mathbb{M}}^{\sf fr}(g,n)$.

\begin{Proposition}~\label{prop:colimit}
There exists a quasi-isomorphism of chain complexes $$\pi^\sharp: {\sf Tot}\big(F\widehat{\mathbb{M}}^{\sf fr}(g,n)_\bullet\big) \ra F\widehat{\mathbb{M}}^{\sf fr}(g,n)$$ depicted as
\[\begin{CD}
F\widehat{\mathbb{M}}^{\sf fr}(g,n)_0 @<\delta<< F\widehat{\mathbb{M}}^{\sf fr}(g,n)_1 @<\delta<< \cdots\\
@V \pi^\sharp VV  @.  @. \\
F\widehat{\mathbb{M}}^{\sf fr}(g,n)  @.  @. 
\end{CD}\]
\end{Proposition}

\begin{proof}
Let us fix a stable graph $G_0\in \Gamma((g,n))_k$  and consider the subspace of the total complex given by
\[ \begin{CD}
F\widehat{\mathbb{M}}^{\sf fr}_{G_0} @<\delta_1+\delta_3<< \displaystyle\varinjlim_{G_0\ra G_1} F\widehat{\mathbb{M}}^{\sf fr}_{G_0\ra G_1} @<\delta_1+\delta_3<<\displaystyle \varinjlim_{G_0\ra G_1\ra G_2} F\widehat{\mathbb{M}}^{\sf fr}_{G_0\ra G_1\ra G_2} @<\delta_1+\delta_3<<  \cdots
\end{CD}\]
which is endowed with the differential $\delta_1+\delta_3$. The natural inclusion map
\[ i: F\widehat{\mathbb{M}}^{\sf fr}(G_0) \hookrightarrow F\widehat{\mathbb{M}}^{\sf fr}_{G_0}\]
which decorates each edge $e$ by the disk cell $D_e$ splits the projection map $$\pi: F\widehat{\mathbb{M}}^{\sf fr}_{G_0}\ra F\widehat{\mathbb{M}}^{\sf fr}(G_0)$$ induced by the projection map $C_*(D) \rightarrow D_e\cdot \mathbb{Q}$ at every edges of $G_0$. We first construct a homotopy operator $$h: \varinjlim_{G_0\ra\cdots\ra G_k} F\widehat{\mathbb{M}}^{\sf fr}_{G_0\ra \cdots \ra G_k} \ra \varinjlim_{G_0\ra\cdots\ra G_{k+1}}F\widehat{\mathbb{M}}^{\sf fr}_{G_0\ra \cdots \ra G_{k+1}}$$ such that it satisfies the deformation retract identity $\id - i\pi = [\delta_1+\delta_3, h]$. To construct such an $h$, consider an element
\[ \alpha \in F\widehat{\mathbb{M}}^{\sf fr}_{G_0\ra\cdots\ra G_k}= \bigotimes_{v\in V_{G_0}}\widehat{\mathbb{M}}_{S^1}^{\sf fr}(g(v),{\sf Leg}(v)) \otimes \bigotimes_{e\in E_{G_0\ra G_k}} C_*(S^1) \otimes \bigotimes_{e'\in E^c_{G_0\ra G_k}} C_*(D).\]
Let $e'_1,\ldots,e'_l$ be edges in $E^c_{G_0\ra G_k}$ such that in $\alpha$ its decoration is in $C_*(S^1)$. Then we set 
\[ h(\alpha) = \alpha \in F\widehat{\mathbb{M}}^{\sf fr}_{G_0\ra\cdots\ra G_k \ra G_{k+1}}\]
where $G_{k+1}$ is obtained from $G_k$ by contracting the edges $e'_1,\ldots,e'_l$. Let us check the identity $[h,\delta_1+\delta_3]=\id-i\pi$. Indeed, if the set $\{ e_1',\ldots,e_l'\}$ is not empty, we have $\alpha= \delta_1 h(\alpha)$ and $\delta_3 h(\alpha)= \sum_{j=1}^k(\alpha)_{G_0\ra\cdots\widehat{G_j}\cdots\ra G_{k+1}}$. While in the other direction, we have
\begin{align*}
h\delta_1(\alpha)&=h\big( (\alpha)_{G_0\ra\cdots\ra G_{k-1}}\big)= (\alpha)_{G_0\ra\cdots\ra G_{k-1}\ra G_{k+1}}\\
h\delta_3(\alpha) &= \sum_{j=1}^{k-1} (\alpha)_{G_0\ra\cdots\widehat{G_j}\cdots\ra G_k \ra G_{k+1}}
\end{align*}
This shows that $[\delta_1+\delta_3,h]=\id$ if $l>0$. Since if $l>0$, we have $i\pi(\alpha)=0$, thus the previous identity is equivalent to the homotopy retract identity. In the case $l=0$, we have $h(\alpha)=0$, and that
\begin{align*}
h\delta_1(\alpha)&=\begin{cases}
h\big( (\alpha)_{G_0\ra\cdots\ra G_{k-1}}\big)= (\alpha)_{G_0\ra\cdots\ra G_{k-1}\ra G_{k}}=\alpha & k\geq 1\\
0 & k=0 
\end{cases}\\
h\delta_3(\alpha) &= 0
\end{align*}
This also implies the desired identity $[h,\delta_1+\delta_3]=\id-i\pi$. 

Finally, we evoke the differential $\partial+\delta_2$ and use homological perturbation formula to compute the induced differential on the Feynman compactification $F\widehat{\mathbb{M}}^{\sf fr}(g,n)$ which matches the definition given in Paragraph~\ref{para:feynman}. Indeed, the perturbed differential turns out to have two terms:
\[ \pi \partial i +  \pi \delta_2 h \partial i\]
The first term corresponds to the boundary map on vertices of stable graphs in the Feynman compactification $F\widehat{\mathbb{M}}^{\sf fr}(g,n) $. The extra differential $\delta$ of $F\widehat{\mathbb{M}}^{\sf fr}(g,n) $ (defined in Equation~\ref{eq:delta}) applied to edges of stable graphs is precisely the second term, i.e.
\[ \delta =  \pi \delta_2 h \partial i.\]
This is a calculation which we illustrate in the case when the underlying stable graph has only one edge $e$ joining two vertices decorated by $\alpha_1$ and $\alpha_2$:
\begin{align*}
\pi\delta_2 h\partial i (D_e\otimes \alpha_1\otimes \alpha_2) &= \pi\delta_2 h\partial (D_e\otimes \alpha_1\otimes \alpha_2) \\
&= \pi \delta_2 (\epsilon_e\otimes \alpha_1\otimes \alpha_2)\\
&= \rho_e(\alpha_1\otimes\alpha_2)
\end{align*}
This agrees with the formula of $\delta$ in Equation~\ref{eq:delta}, which finishes the proof.
\end{proof}

\paragraph We write down the perturbed inclusion and projection maps which will be used later in the paper.
\begin{align}~\label{eq:in-pro}
\begin{split}
i^{\sharp} & = \sum_{k\geq 0} (h\partial)^k i \\
\pi^{\sharp} &= \pi +   \pi \delta_2 h
\end{split}
\end{align}

Completely parallel to the construction of $F\widehat{\mathbb{M}}^{\sf fr}(g,n)_\bullet$, one may also define a simplicial resolution $F\mathbb{M}^{\sf fr}(g,n)_\bullet$ for the Feynman compactification $F\mathbb{M}^{\sf fr}(g,n)$. We also have the quasi-isomorphisms:
\begin{align*}
i^\sharp & : F\mathbb{M}^{\sf fr}(g,n) \ra  {\sf Tot}\big( F{\mathbb{M}}^{\sf fr}(g,n)_\bullet \big)\\
\pi^\sharp & : {\sf Tot}\big( F{\mathbb{M}}^{\sf fr}(g,n)_\bullet \big) \ra F\mathbb{M}^{\sf fr}(g,n)
\end{align*}

\section{Comparison with the Deligne-Mumford operad}~\label{sec:isomorphism}

In this section, we prove the following theorem.

\begin{Theorem}~\label{thm:dm-operad}
There is an isomorphism of modular operads $$H_*(F\mathbb{M}^{\sf fr}) \cong H_*(\overline{{M}},\mathbb{Q})$$ where $\overline{{M}}$ is the Deligne-Mumford modular operad~\cite[Section 6.2]{GetKap}. 
\end{Theorem}

\medskip
\medskip
The proof of this theorem occupies the rest of the section. We briefly outline the main idea here. By results of the previous section, we have isomorphisms:
\begin{align*}
H_*\big( F\mathbb{M}^{\sf fr}(g,n)\big) &\cong H_*\Big({\sf Tot}\big(  F\mathbb{M}^{\sf fr}(g,n)_\bullet\big)\Big)\\
H_*\big( F\widehat{\mathbb{M}}^{\sf fr}(g,n)\big) &\cong H_*\Big({\sf Tot}\big( F\widehat{\mathbb{M}}^{\sf fr}(g,n)_\bullet\big) \Big) \cong H_*(\overline{M}_{g,n})
\end{align*}
Thus in order to relate the first line with the second line, we shall construct a quasi-isomorphism
\[ \mathbb{I}: {\sf Tot}\big(  F\mathbb{M}^{\sf fr}(g,n)_\bullet\big)\cong {\sf Tot}\big( F\widehat{\mathbb{M}}^{\sf fr}(g,n)_\bullet\big) \]
Combining with the previous isomorphisms yields the isomorphism in Theorem~\ref{thm:dm-operad}. One then checks the compatibility of the modular operad structures.

\paragraph{{\bf Comparison between $F\widehat{\mathbb{M}}^{\sf fr}(g,n)_\bullet$ and $F\mathbb{M}^{\sf fr}(g,n)_\bullet$.}} As was pointed in Paragraph~\ref{para:ex2}, the inclusion map
\[I_{g,n} : \mathbb{M}^{\sf fr}(g,n)\hookrightarrow \widehat{\mathbb{M}}^{\sf fr}(g,n)\]
does not respect the operadic composition map. This implies that the induced inclusion map between $F\widehat{\mathbb{M}}^{\sf fr}(g,n)_\bullet$ and $F\mathbb{M}^{\sf fr}(g,n)_\bullet$ does not respect the simplicial structure. More precisely, the simplicial boundary map $\delta_2$ is not compatible with the inclusion map. However, the failure of this compatibility is actually homotopically trivial. Consider the simplest situation with $h: G \ra \star_{g,n}$ with $G$ a stable graph with two vertices $v_1$, $v_2$ and an unique edge $e$. Consider the following (non-commutative) diagram:
\[\begin{CD}
M_{g(v_1),\Leg(v_1)}^{\sf fr}\times M_{g(v_2),\Leg(v_2)}^{\sf fr} @> \eta_h >> M_{g,n}^{\sf fr}\\
@V I VV       @VV I V\\
\widehat{M}_{g(v_1),\Leg(v_1)}^{\sf fr}\times \widehat{M}_{g(v_2),\Leg(v_2)}^{\sf fr} @> \xi_h >> \widehat{M}_{g,n}^{\sf fr}
\end{CD}\]
The above diagram is commutative up to homotopy. Indeed, let $\phi_{1}: D_z\ra U(x)$ and $\phi_{2}: D_w\ra U(y)$ be the two framings to be sewed. Fix $0<t\leq 1$, we may scale the two framings by setting
\[ \phi_1^t(z):= \phi_1(t\cdot z),\;\;\; \phi_2^t(w):=\phi_2(t\cdot w).\]
The scaled framings are given by $\phi_1^t: D_z \ra \phi_1(t\cdot D)\subset U(x)$ and $ \phi_2^t: D_w \ra \phi_1(t\cdot D)\subset U(x)$. Using the scaled framings we obtain a $t$-dependent sewing map
\[ \eta_h^t: M_{g(v_1),\Leg(v_1)}^{\sf fr}\times M_{g(v_2),\Leg(v_2)}^{\sf fr} \ra M_{g,n}^{\sf fr}\]
It is clear that $\eta_h^1=\eta_h$.

\begin{Lemma}~\label{lem:limit-sewing}
Let the notations be as above. Then we have 
\[ \lim_{t\ra 0} I\circ \eta_h^t = \xi_h\circ I.\]
\end{Lemma}

\begin{proof}
It is clear that on the underlying stable curve, as $t\ra 0$ the limiting curve becomes a nodal curve with the node given by the image of $x$ and $y$. The limiting decoration at this node is then given by the unit tangent direction vector of
\[ d\phi_{1}^t(\frac{d}{dz})\otimes d\phi_{2}^t(\frac{d}{dw})=t^2 d\phi_{1}(\frac{d}{dz})\otimes d\phi_{2}(\frac{d}{dw}).\]
Thus the unit tangent direction is independent of $t$, and thus have a limiting vector given by the unit tangent direction vector of $d\phi_{1}(\frac{d}{dz})\otimes d\phi_{2}(\frac{d}{dw})$, which is precisely the definition of $\xi_h\circ I$.
\end{proof}

Using the geometric homotopy $\eta_h^t$, we obtain a homotopy operator of homological degree one:
\[ \mathbb{I}_1 :  \mathbb{M}^{\sf fr}(G) \ra \widehat{\mathbb{M}}^{\sf fr}(\star_{g,n}).\]
The lemma above implies that $[\partial,\mathbb{I}_1]=I\circ \eta_h - \xi_h\circ I$. 

In general, for each $k$-step contraction, consider the following diagram:
\[\begin{CD}
F{\mathbb{M}}^{\sf fr}_{G_0\ra\cdots\ra G_k} @>\delta_2>> F{\mathbb{M}}^{\sf fr}_{G_1\ra\cdots\ra G_k}\\
@V I VV     @V I VV\\
F\widehat{\mathbb{M}}^{\sf fr}_{G_0\ra\cdots\ra G_k} @>\delta_2>> F\widehat{\mathbb{M}}^{\sf fr}_{G_1\ra\cdots\ra G_k}
\end{CD}\] 
We may apply the previous construction to all the edges in $E_{G_0\ra G_1}$ with the same parameter $t\in [0,1]$ to yield a homotopy operator still denoted by
\[ \mathbb{I}_1: F{\mathbb{M}}^{\sf fr}_{G_0\ra\cdots\ra G_k} \ra F\widehat{\mathbb{M}}^{\sf fr}_{G_1\ra\cdots\ra G_k}\] 
of homological degree one such that $[\partial,\mathbb{I}_1]=I\circ \delta_2 - \delta_2 \circ I$. 
Using these homotopy operators and its higher extensions (introduced in the proof below), we prove the following

\begin{Theorem}~\label{thm:bbI}
There exists a quasi-isomorphism of chain complexes
\[ \mathbb{I}: {\sf Tot}\big( F{\mathbb{M}}^{\sf fr}(g,n)_\bullet\big)\ra  {\sf Tot}\big( F\widehat{\mathbb{M}}^{\sf fr}(g,n)_\bullet\big)\]
\end{Theorem}

\begin{proof}
We inductively construct a morphism of chain complexes of the form
\[ \mathbb{I}=\sum_{d\geq 0} \mathbb{I}_d : {\sf Tot}\big( F{\mathbb{M}}^{\sf fr}(g,n)_\bullet\big)\cong {\sf Tot}\big( F\widehat{\mathbb{M}}^{\sf fr}(g,n)_\bullet\big)\]
with $\mathbb{I}_0:=I$. Recall that the total differential on both total complexes is of the form $\partial + \delta$ with $\partial$ the internal differential, and $\delta=\delta_1+\delta_2+\delta_3$ the simplicial boundary maps. From the definition, it is easy to verify that $[\delta_1,\mathbb{I}_0]=[\delta_3,\mathbb{I}_0]=0$. But the boundary map $\delta_2$ is not compatible with $\mathbb{I}_0$. Its commutator is precisely bounded by the homotopy $\mathbb{I}_1$ defined in the previous paragraph.

In general, for each $d\geq 1$, we shall define a map of homological degree $d$ of the form
\[ \mathbb{I}_d: (F{\mathbb{M}}^{\sf fr})_{G_0\ra\cdots\ra G_k} \ra (F\widehat{\mathbb{M}}^{\sf fr})_{G_d\ra\cdots\ra G_k}\]
such that 
\begin{align}~\label{eq:I_d}
\begin{split}
[\partial,\mathbb{I}_{d+1}] & =-[\delta_2+\delta_3,\mathbb{I}_d],\\
[\delta_1,\mathbb{I}_{d}] & = 0
\end{split}
\end{align}
Indeed, associated with a $d$-step stable graph contraction $G_0\ra\cdots\ra G_d$ we have a family of sewing maps depending on $d$-variables $t_1,\cdots,t_d$
\[ \eta_{G_0\ra\cdots\ra G_d}^{t_1,\cdots,t_d}: F{\mathbb{M}}^{\sf fr}(G_0)\ra F\widehat{\mathbb{M}}^{\sf fr}(G_d), \;\; 0\leq t_1\leq\ldots\leq t_d\leq 1,\]
where we use the $t_i$-variable to scale the sewing operation at an edge that got contracted by $G_{i-1}\ra G_i$. See Figure $1$ for an illustration of these higher homotopy maps.

\begin{figure}~\label{fig:homotopy}
\begin{center}
\includegraphics[scale=.4]{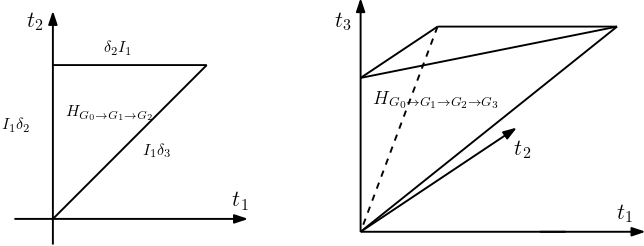}
\end{center}
\caption{Higher homotopies}
\end{figure}

We use this $d$-parameter family to induce a higher homotopy operator:
\[ \mathbb{I}_d: (F{\mathbb{M}}^{\sf fr})_{G_0\ra\cdots\ra G_k} \ra (F\widehat{\mathbb{M}}^{\sf fr})_{G_d\ra \cdots\ra G_k}.\]
To see the above identities~\eqref{eq:I_d} hold, observe that the facets defined by $t_1=0$ and $t_d=1$ corresponds to $-[\delta_2,\mathbb{I}_{d-1}]$, the intermediate facets defined by $t_i=t_{i+1}, \; 1\leq i \leq d-1$ are $\mathbb{I}_{d-1}\delta_3=-[\delta_3,\mathbb{I}_{d-1}]$. Here we used that $\delta_3\mathbb{I}_{d-1}=0$ since the image of $\mathbb{I}_{d-1}$ lies in the index $G_{d-1}\ra G_{d}$ while $\delta_3$ vanishes on such type of component. Putting the maps $\{\mathbb{I}_d\}_{d\geq 0}$ together, we obtain a map of chain complexes
\[ \mathbb{I}:=\sum_{d\geq 0} \mathbb{I}_d: {\sf Tot}\big( F{\mathbb{M}}^{\sf fr}(g,n)_\bullet\big)\ra {\sf Tot}\big( F\widehat{\mathbb{M}}^{\sf fr}(g,n)_\bullet\big).\]
To argue that $\mathbb{I}$ is a quasi-isomorphism, we consider the increasing filtration on both sides defined by that the $k$-th filtration consists of elements of simplicial degree less than or equal to $k$.  Observe that $\mathbb{I}$ preserves the filtration. Thus it induces a morphism of the associated spectral sequences on the two total complexes. But we already have an isomorphism in the $E^1$-page since the map $\mathbb{I}_0$ is induced by the natural inclusion maps $M_{g,n}^{\sf fr} \hookrightarrow \widehat{M}_{g,n}^{\sf fr}$ which are homotopy equivalences.
\end{proof}

\paragraph{{\bf Operadic compositions.}}~\label{para:composition} To summarize the previous discussions, we have the following diagram
\[\begin{CD}
{\sf Tot}\big( F{\mathbb{M}}^{\sf fr}(g,n)_\bullet\big) @>\mathbb{I}>>  {\sf Tot}\big( F\widehat{\mathbb{M}}^{\sf fr}(g,n)_\bullet\big) \\
@ A i^\sharp AA    @V \pi^\sharp VV \\
F{\mathbb{M}}^{\sf fr}(g,n)           @.   F\widehat{\mathbb{M}}^{\sf fr}(g,n) 
\end{CD}\]
where all the arrows are quasi-isomorphisms. To prove analyze the operadic compositions in Theorem~\ref{thm:dm-operad}, we shall define a canonical projection map
\[ p: {\sf Tot}\big( F\widehat{\mathbb{M}}^{\sf fr}(g,n)_\bullet\big) \ra C_*(\overline{M}_{g,n}).\] This map is only non-zero at simplicial degree zero, in which case we need to define a map
\[ p: (F\widehat{\mathbb{M}}^{\sf fr})_G \ra C_*(\overline{M}_{g,n})\]
for each stable graph $G\in \Gamma((g,n))$. Let $\alpha\in (F\widehat{\mathbb{M}}^{\sf fr})_G $ with its underlying stable graph $G$. Explicitly, such an element is of the form
\[  \alpha= \prod_{e\in E_G} \alpha_e  \otimes \prod_{v\in V_G} \alpha_v \in \prod_{e\in E_G} C_*(D_e) \otimes \prod_{v\in V_G} C_*(\widehat{M}_{g(v),n(v)}^{\sf fr})[u_1^{-1},\ldots,u_{n(v)}^{-1}]\]
At a vertex $v\in V_G$, its vertex decoration is of the form
\[ \alpha_v= \alpha_v^{[0]} + O(u^{-1}),\]
with $\alpha_v^{[0]}\in C_*(M_{g(v),n(v)}^{\sf fr})$ the constant term in $\alpha_v\in C_*(\widehat{M}_{g(v),n(v)}^{\sf fr})[u_1^{-1},\ldots,u_{n(v)}^{-1}]$. Let us denote by
\[ p_v: C_*(\widehat{M}_{g(v),n(v)}^{\sf fr})[u_1^{-1},\ldots,u_{n(v)}^{-1}] \ra C_*(\overline{M}_{g(v),n(v)})\]
the map obtained by first projection onto the constant term, followed by the forgetful map ${\sf Forget}: \widehat{M}_{g(v),n(v)}^{\sf fr} \ra \overline{M}_{g(v),n(v)}$ which forgets both the framings and the nodal decorations. For each edge $e\in E_G$, denote by
\[ p_e: C_*(D_e) \ra \mathbb{Q}\]
the projection map onto $C_0(D)=\mathbb{Q}$. Then the map $p$ is defined by
\[ p(\alpha)= \xi_G\big( \prod_{e\in E_G} p_e(\alpha_e)  \otimes \prod_{v\in V_G} p_v(\alpha_v) \big)\]
where $\xi_G$ is the sewing map in Equation~\eqref{eq:xi}. Using the map $p$, we consider the composition
\[  p \mathbb{I} i^\sharp: F{\mathbb{M}}^{\sf fr}(g,n) \ra C_*( \overline{M}_{g,n}).\]

\begin{Proposition}
The induced map on homology
\[ p \mathbb{I} i^\sharp: H_*( F{\mathbb{M}}^{\sf fr} ) \ra H_*( \overline{M})\]
is a morphism of modular operads.
\end{Proposition}

\begin{proof}
We shall only deal with the first type composition corresponding to a loop contraction $c_{ij}$. The case of non-loop contraction is analogous. In the loop case, we want to prove that the following diagram is commutative up to homotopy:
\[\begin{CD}
F{\mathbb{M}}^{\sf fr}(g,n+2)  @>F{\mathbb{M}}^{\sf fr}(c_{ij})>> F{\mathbb{M}}^{\sf fr}(g+1,n)\\
@V p\mathbb{I} i^\sharp VV      @V p\mathbb{I} i^\sharp VV\\
C_*(\overline{M}_{g,n+2}) @>C_*(\overline{M})(c_{ij})>> C_*(\overline{M}_{g+1,n})
\end{CD}\]
Let $\alpha\in F{\mathbb{M}}^{\sf fr}(g,n+2)$ be an element with its underlying stable graph $G\in \Gamma((g,n+2))$. Let us explicitly write down the composition
 \[ p \mathbb{I} i^\sharp:  F\mathbb{M}^{\sf fr}(g,n)\ra C_*(\overline{M}_{g,n})\]
using Equations of $p$, $\mathbb{I}$ and Equation~\eqref{eq:in-pro}. We have
\begin{align*}
\mathbb{I}i^\sharp (\alpha) & = \sum_{k\geq 0} \sum_{(e_1,\ldots,e_k)} \mathbb{I}_k \big( \prod_{j=1}^k \epsilon_{e_j} \otimes \prod_{e\in E_G-\{e_1,\ldots,e_k\}} D_e \otimes \prod_{v\in V_G} \alpha_v\big)_{G\ra G/e_1 \ra \cdots \ra G/(e_1,\ldots,e_k)}\\
\end{align*}
Since $p_e$ vanishes on $D_e$, we have
\[ p\mathbb{I}i^\sharp (\alpha) = \sum_{(e_1,\ldots,e_K)} \mathbb{I}_K \big( \prod_{e\in E_G} \epsilon_{e} \otimes \prod_{v\in V_G} \alpha_v\big)_{G\ra G/e_1 \ra \cdots \ra G/(e_1,\ldots,e_K)}\]
where the summation is over orderings of the edge set $E_G$, and $K=|E_G|$. Observe that in the definition of $\mathbb{I}_K$, we required the sewing parameters $t_1,\ldots,t_K$ be ordered as $0\leq t_1\leq \ldots,\leq t_K\leq 1$. Summing over ordering of the edges  gets rid of this ordering condition. Thus, the above expression can be rewritten as
\begin{equation}~\label{eq:polydisk}
p\mathbb{I}i^\sharp (\alpha) = {\sf Forget}\Big( \eta_G\big(\underbrace{[0,1]\times\cdots\times[0,1]}_{\mbox{$K$ copies}} \times ( \prod_{e\in E_G} B_e) (\prod_{v\in V_G} \alpha_{g(v),n(v)}^{[0]}) \big)\Big)
\end{equation}
with
\[ \eta_G: \underbrace{[0,1]\times\cdots\times[0,1]}_{\mbox{$K$ copies}} \times \prod_{v\in V_G} M^{\sf fr}_{g(v),n(v)} \ra \widehat{M}_{g,n}^{\sf fr}\]
the sewing map depending on $K$ parameters as in the proof of Theorem~\ref{thm:bbI}.
Intuitively speaking, for each $e\in E_G$, we act on the chain $\prod_{v\in V_G} \alpha_{g(v),n(v)}^{[0]}$ by the circle action $B_e$ and also sewing with the radius parameter $t_e\in [0,1]$. The forgetful map precisely collapses the circle at $t_e=0$, which forms a disk $[0,1]\times S^1 / {0}\times S^1$ over each edge $e\in E_G$. 

To this end, we shall finish the proof using the explicit formula of $p\mathbb{I} i^\sharp$ above. Indeed, let $\alpha\in F{\mathbb{M}}^{\sf fr}(g,n+2)$ be such that at the half-edges labeled by $i$ and $j$ the $u$-powers are $u^{-k_i}$ and $u^{-k_j}$. If either $k_i$ or $k_j$ is strictly positive, the resulting composition is trivial in homology. In the interesting case when $k_i=k_j=0$, denote by $G_{ij}\in \Gamma(g+1,n)$ the graph obtained from $G$ by sewing the two leaves indexed by $i$ and $j$. Then we see that the two compositions $p\mathbb{I} i^\sharp F\mathbb{M}^{\sf fr} (c_{ij})(\alpha)$ and $C_*(\overline{M})(c_{ij})p\mathbb{I} i^\sharp(\alpha)$ are related by the homotopy
\[  \eta_{G_{ij}}: \underbrace{[0,1]\times\cdots\times[0,1]}_{\mbox{$K+1$ copies}} \times \prod_{v\in V_G} M^{\sf fr}_{g(v),n(v)} \ra \widehat{M}_{g,n}^{\sf fr}\]
This finishes the proof of proposition, which combined with Theorem~\ref{thm:bbI} yields Theorem~\ref{thm:dm-operad}.
\end{proof}

\section{String vertices and fundamental classes}~\label{sec:fundamental}

In this section, we use the Feynman compactification construction to provide a geometric extension of the algebraic trivialization of circle actions in~\cite{CT}. Roughly speaking, we construct a compactified version denoted by $\overline{\mathfrak{g}}$ 
of Sen-Zwiebach's DGLA $\mathfrak{g}$ (reviewed in the beginning of the section) so that we have an inclusion of DGLA's
\[ j: \mathfrak{g}\ra \overline{\mathfrak{g}}.\]
Then we construct an $L_\infty$ quasi-isomorphism $\mathcal{K}: \overline{\mathfrak{g}} \ra \overline{\mathfrak{g}}^{\sf triv}$ where the DGLA $\overline{\mathfrak{g}}^{\sf triv}$ has the same underlying space as $\overline{\mathfrak{g}}$, but its BV differential and the Lie bracket are both zero. We prove the identity~\eqref{eq:fundamental} stated in the introduction that the push-forward of the string vertex $\mathcal{K}_* j_* \mathcal{V}$ consists of the fundamental class of $[\overline{M}_{g,n}/S_n]$. 

\paragraph{{\bf DGLA's.}}~\label{para:DGLA} Recall the construction of Sen-Zwiebach's DGLA:
\[     \mathfrak{g} := \big(\bigoplus_{(g,n)} \mathbb{M}^{\sf fr}(g,n)_{S_n\ltimes (hS^1)^n} [1]\big)[[\hbar,\lambda]]\]
with two formal variables $\hbar$ and $\lambda$ both of homological degree $-2$. By construction, its differential is $\partial+uB+\hbar \Delta$ with $\Delta$ defined using the twisted self-sewing map~\eqref{eq:self-twisted-sewing}, and its Lie bracket  is $\{-,-\}$ defined using the twisted sewing map~\eqref{eq:twisted-sewing}. We refer to~\cite[Section 3]{CCT} for more details of this construction.

Similar to the definition of $\mathfrak{g}$, one may define a ``compactified" version of $\mathfrak{g}$ by putting
\[  \overline{\mathfrak{g}} :=  \big(\bigoplus_{(g,n)} F\mathbb{M}^{\sf fr}(g,n)_{S_n} [1]\big)[[\hbar,\lambda]]\]
with the differential $\partial+uB+\delta+\hbar \Delta$ (with $\delta$ as in Equation~\eqref{eq:delta}) and the Lie bracket $\{-,-\}$ defined in the same way as $\mathfrak{g}$. Observe that there is an inclusion map
\[ j: \mathfrak{g}\ra \overline{\mathfrak{g}}\]
of DGLA's onto the $\star_{g,n}$-component of each $F\mathbb{M}^{\sf fr}(g,n)_{S_n}$.

We also define a trivialized version of $\overline{\mathfrak{g}}$, denoted by $\overline{\mathfrak{g}}^{\sf triv}$ which has the same underlying graded vector space as $\overline{\mathfrak{g}}$, but is endowed with differential $\partial+uB+\delta$ and zero Lie bracket.

\begin{Theorem}
There exists an $L_\infty$ quasi-isomorphism of DGLA's:
\[ \mathcal{K}: \overline{\mathfrak{g}} \ra \overline{\mathfrak{g}}^{\sf triv}\]
\end{Theorem}

This theorem is the geometric ``lift" of the algebraic version proved in~\cite[Theorem 4.2]{CT}. The proof is in complete parallel as well. In the following, we describe the construction of $\mathcal{K}$. Roughly speaking, the map is simply inserting elements of $\overline{\mathfrak{g}}$ on vertices of stable graphs, as illustrated in Figure $2$.

\paragraph{{\bf Construction of $\mathcal{K}$.}} We proceed to construct an $L_\infty$-morphism $\mathcal{K}$. First, we introduce more notations about graphs. Denote by $\Gamma(g,n)_m$ (respectively $\Gamma((g,n))_m$) the set of labeled (respectively stable) graphs with $m$ vertices. For a labeled graph $G\in \Gamma(g,n)_m$, a marking of $G$ is a bijection $$f: \{1,\cdots,m\} \ra V_G.$$ An isomorphism between two marked and labeled graphs is an isomorphism of the underlying labeled graphs that also preserves the marking map. Denote by $\widetilde{\Gamma(g,n)}_m$ (respectively $\widetilde{\Gamma((g,n))}_m$) the set of isomorphism classes of marked (respectively stable) graphs. 

For each integer $m\geq 1$, we shall define a degree zero linear map
\[ \mathcal{K}_m: \sym^m (\overline{\mathfrak{g}}[1]) \ra \overline{\mathfrak{g}}^{\sf triv}[1]\]
The shift by one of $ \overline{\mathfrak{g}}$ is given by
$$ \overline{\mathfrak{g}}[1] =  \big(\bigoplus_{(g,n)} F\mathbb{M}^{\sf fr}(g,n)_{S_n} [2]\big)[[\hbar,\lambda]]$$
For each marked labeled graph $(G,f)\in \widetilde{\Gamma(g,n)}_m$ with $f: \{ 1,\ldots, m\} \ra V_G$ a marking on the vertices of $G$,  denote by $v_j:=f(j)$ and $n_j:=|{\sf Leg}(v_j)|$. Define a $\mathbb{Q}$-linear map 
\[ \mathcal{K}_{(G,f)}:\bigotimes_{j=1}^m F\mathbb{M}^{\sf fr}(g(v_j),n_j)_{S_{n_j}}[2]\cdot \hbar^{g(v_j)} \ra F\mathbb{M}^{\sf fr}(g,n)_{S_{n}}[2]\cdot \hbar^g.\]

\begin{figure}~\label{fig:K-construction}
\begin{center}
\includegraphics[scale=.6]{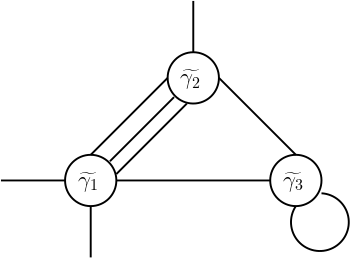}
\end{center}
\caption{Illustration of $\mathcal{K}_{(G,f)}(\gamma_1,\gamma_2,\gamma_3)$}
\end{figure}

For each $1\leq j\leq m$, let
\[ \gamma_j \in \big( D(H_j)\bigotimes \bigotimes_{v\in V_{H_j}} \mathbb{M}_{S^1}^{\sf fr}(g(v),{\sf Leg}(v))\big)_{{\sf Aut}(H_j)}\]
be an element of $F\mathbb{M}^{\sf fr}(g(v_j),n_j)$, with its underlying stable graph given by $H_j\in \Gamma(g(v_j),n_j)$. In order to insert the graphs $H_j$'s into $G$, we need to make identification of leaves of $H_j$ with ${\sf Leg}(v_j)$. There is no canonical choice for such an identification, thus we take the symmetrization of the $\gamma_j$'s. We set
\[ \gamma_j^{\sym} := \sum_{\sigma_j: \{1,\ldots,n_j\} \ra {\sf Leg}(v_j)}  \gamma_j^{\sigma_j},\]
where the notation $\gamma_j^{\sigma_j}$ is the same element as $\gamma_j$, but  labeled by an isomorphism $\sigma_j: \{1,\ldots,n_j\} \ra {\sf Leg}(v_j)$. Fixing the identifications $\sigma_1,\ldots,\sigma_m$ on vertices of $G$, denote by
\[ G(H_1^{\sigma},\ldots,H_m^{\sigma_m})\]
the stable graph obtained by inserting $H_j$ at vertex $v_j$ through the isomorphism $\sigma_j$. With this preparation, we define
\[  \mathcal{K}_{(G,f)}(\gamma_1,\ldots,\gamma_m):= \sum_{\sigma_1,\ldots,\sigma_m} D(G)\otimes \gamma_1^{\sigma_1}\otimes\cdots\otimes \gamma_m^{\sigma_m}.\]
This makes sense because
\begin{itemize}
\item There is a natural isomorphism
\[ D(G(H_1^{\sigma_1},\ldots,H_m^{\sigma_m}))\cong D(G)\otimes D(H_1)\otimes\cdots\otimes D(H_m)\]
\item We have $V_{G(H_1^{\sigma_1},\ldots,H_m^{\sigma_m})}= \coprod_{j=1}^m V_{H_j}$.
\item This map descends to coinvariants by graph automorphisms. Indeed, given an automorphism $\rho_j: H_j \ra H_j$, it induces an isomorphism
\[ G(H_1^{\sigma_1},\ldots,H_m^{\sigma_m}) \ra G(H_1^{\sigma_1},\ldots,H_j^{\rho_j. \sigma_j},\ldots,H_m^{\sigma_m}),\]
where $\rho_j.\sigma_j$ is the identification given by first applying $\sigma_j$ then followed by the automorphism $\rho_j$ acting on leaves of $H_j$.
\item This map is of homological degree zero after the shift by $[2]$. To see this, observe that at vertices the degrees are simply that in $\gamma_j$'s, but since $D(G(H_1^{\sigma_1},\ldots,H_m^{\sigma_m}))\cong D(G)\otimes D(H_1)\otimes\cdots\otimes D(H_m)$ the term $D(G)$ has degree $2|E_G|$. Putting this together with the shift by $[2]$ and the powers of $\hbar$, we have
\begin{align*}
 \deg(\mathcal{K}_{(G,f)})& = 2|E_G|-2|V_G|+2-2(g-g(v_1)-\cdots - g(v_m))\\
 &= 2(-1+h_1(G))+2-2(g-g(v_1)-\cdots - g(v_m))\\
 &= 2(h_1(G)+g(v_1)+\cdots + g(v_m)-g)=0
 \end{align*}
\end{itemize}
Finally, extending $\mathcal{K}_{(G,f)}$ by $\lambda$-linearity, we define for each $m\geq 1$ a linear (in both $\lambda$ and $\hbar$ variables) by setting
\begin{equation}~\label{eq:K}
 \mathcal{K}_m:= \sum_{\substack{(g,n)\\(G,f)\in\widetilde{\Gamma(g,n)}_m}} \frac{1}{|{\sf Aut}(G,f)|}\cdot \mathcal{K}_{(G,f)}
\end{equation}

\paragraph{{\bf Proof of Equation~\eqref{eq:fundamental}.}} In this subsection, we prove the main identity~\eqref{eq:fundamental} that the push-forward of the String vertex $\mathcal{V}= \sum_{(g,n)} \mathcal{V}_{g,n} \hbar^g\lambda^{2g-2+n}$ via the composition $L_\infty$ map
\[ \mathfrak{g} \hookrightarrow \overline{\mathfrak{g}} \stackrel{\mathcal{K}}{\longrightarrow} \overline{\mathfrak{g}}^{\sf triv}\]
gives exactly the fundamental classes $ [\overline{M}_{g,n}/S_n]$.

\begin{Theorem}~\label{thm:main}
Under the isomorphism in Theorem~\ref{thm:dm-operad}, we have the following formula expressing the fundamental class of $\overline{M}_{g,n}/S_n$ in terms of String vertices:
\[ [\overline{M}_{g,n}/S_n]= (\mathcal{K}_*\mathcal{V})_{g,n}= \sum_{G\in \Gamma((g,n))} \frac{1}{|{\sf Aut}(G)|} \prod_{e\in E_G} D_e \otimes \prod_{v\in V_G} \mathcal{V}^{\sym}_{g(v),|{\sf Leg}(v)|}.\]
\end{Theorem}

\begin{proof}
Define a filtration on the chain complex $\widehat{\mathbb{M}}^{\sf fr}(g,n)$ by the setting $\partial_k\widehat{\mathbb{M}}^{\sf fr}(g,n)$ to be chains that are supported on the locus such that the underlying stable Riemann surfaces have at least $k$ nodes. Similarly, we also have the nodal filtration on the homotopy quotient $\partial_k\widehat{\mathbb{M}}_{S^1}^{\sf fr}(g,n) \; (k\geq 0)$, and further $S_n$-quotient $\partial_k\widehat{\mathbb{M}}_{S^1}^{\sf fr}(g,n)_{S_n}$. Following Costello~\cite[Theorem 9.0.7]{Cos2}, one may form a DGLA on the homology of the associated graded
\[ \mathfrak{n}:= \bigoplus_{(g,n),k\geq 0} H_*\big( {\sf Gr}^k_{\partial_\bullet}\widehat{\mathbb{M}}_{S^1}^{\sf fr}(g,n)_{S_n}\big)[1] [[\hbar,\lambda]]\]
Its differential and Lie bracket are denoted by $\partial + \hbar\Delta$ and $\{-,-\}$ respectively. The boundary map $\partial:H_d\big( {\sf Gr}^k_{\partial_\bullet}\widehat{\mathbb{M}}_{S^1}^{\sf fr}(g,n)_{S_n}\big) \ra H_{d-1}\big( {\sf Gr}^{k+1}_{\partial_\bullet}\widehat{\mathbb{M}}_{S^1}^{\sf fr}(g,n)_{S_n}\big)$ is the boundary map induced by the sequence
\[ 0\ra  {\sf Gr}^{k+1}_{\partial_\bullet}\widehat{\mathbb{M}}_{S^1}^{\sf fr}(g,n)_{S_n} \ra \partial_k\widehat{\mathbb{M}}_{S^1}^{\sf fr}(g,n)_{S_n}/\partial_{k+2}\widehat{\mathbb{M}}_{S^1}^{\sf fr}(g,n)_{S_n}\ra   {\sf Gr}^{k}_{\partial_\bullet}\widehat{\mathbb{M}}_{S^1}^{\sf fr}(g,n)_{S_n} \ra 0\]
The two operations $\Delta$ and $\{-,-\}$  are both defined using the twisted sewing operations. Observe that when $k=0$ we have
\begin{align*}
& H_*\big( {\sf Gr}^0_{\partial_\bullet}\widehat{\mathbb{M}}_{S^1}^{\sf fr}(g,n)_{S_n}\big)\\
= & H_*\big( \widehat{\mathbb{M}}_{S^1}^{\sf fr}(g,n)_{S_n}, \partial_1\widehat{\mathbb{M}}_{S^1}^{\sf fr}(g,n)_{S_n}\big)\\
\cong & H_*\big( \overline{M}_{g,n}/S_n, \partial_1 \overline{M}_{g,n}/S_n\big)
\end{align*}
The last isomorphism is induced by the canonical map $\widehat{M}_{g,n} \ra \overline{M}_{g,n}$ that forgets the decoration at all nodes. The orbifold fundamental class $[\overline{M}_{g,n}/S_n]$, under this isomorphism, corresponds to the fundamental class $[\widehat{M}_{g,n}/S_n]\in H_{6g-6+2n}(\widehat{M}_{g,n}/S_n, \partial_1\widehat{M}_{g,n}/S_n)$. It was shown in {\sl Loc. Cit.} that the element
\[ \mathcal{X}:= \sum_{(g,n)} [\widehat{M}_{g,n}/S_n] \hbar^g \lambda^{2g-2+n}\]
of the DGLA $\mathfrak{n}$ satisfies the Maurer-Cartan equation and is uniquely determined by the equation after fixing the initial condition at $(g,n)=(0,3)$. 

To this end, let us consider the following commutative diagram
\[ \begin{CD}
F{\mathbb{M}}^{\sf fr}(g,n) @>i^\sharp>> {\sf Tot}\big( F{\mathbb{M}}^{\sf fr}(g,n)_\bullet\big) @>\mathbb{I}>>  {\sf Tot}\big( F\widehat{\mathbb{M}}^{\sf fr}(g,n)_\bullet\big) @>p>> C_*(\overline{M}_{g,n})\\
@. @. @V{\sf pr} VV @A{\sf Forget} AA\\
@. @.  \widehat{\mathbb{M}}^{\sf fr}_{S^1}(g,n)@>q >> C_*(\widehat{M}_{g,n})
\end{CD}\]
where the downward arrow is the projection map onto the $\star_{g,n}$ component, the upward arrow is the map that forgets the nodal decoration, and the bottom horizontal arrow $q$ is the canonical quotient map followed by forgetting  framings at punctures. Recall from Paragraph~\ref{para:composition} the isomorphism in Theorem~\ref{thm:dm-operad} is induced by the chain map $p\mathbb{I} i^\sharp$. By the previous discussion, it suffices to prove that
\[q {\sf pr} \mathbb{I} i^\sharp\big(\mathcal{K}_*\mathcal{V}\big)_{g,n}=q {\sf pr} \mathbb{I} i^\sharp \big( \sum_{G\in \Gamma((g,n))} \frac{1}{|{\sf Aut}(G)|} \prod_{e\in E_G} D_e \otimes \prod_{v\in V_G} \mathcal{V}^{\sym}_{g(v),|{\sf Leg}(v)|}\big)\]
satisfies the Maurer-Cartan equation in $\mathfrak{n}$. Using the identity~\eqref{eq:polydisk} we have
\begin{align*}
& q{\sf pr} \mathbb{I}i^\sharp \big(\mathcal{K}_*\mathcal{V}\big)_{g,n}\\
= & q\Big( \sum_{G\in \Gamma((g,n))}\frac{1}{|{\sf Aut}(G)|} \cdot \eta_G\big(\underbrace{[0,1]\times\cdots\times[0,1]}_{\mbox{$|E_G|$ copies}} \times ( \prod_{e\in E_G} B_e) (\prod_{v\in V_G}  \mathcal{V}^{\sym,[0]}_{g(v),n(v)}) \big)\Big)
\end{align*}
To finish the proof, we show that the Maurer-Cartan equation satisfied by $\mathcal{V}_{g,n}$'s implies the desired Maurer-Cartan equation of $q{\sf pr} \mathbb{I}i^\sharp \big(\mathcal{K}_*\mathcal{V}\big)_{g,n}$'s. Indeed, the leading coefficients $\mathcal{V}_{g,n}^{[0]}$ of string vertices satisfies the equation:
\begin{equation}~\label{eq:leading-mc}
 \partial \mathcal{V}_{g,n}^{[0]} = -\frac{1}{2} \sum_{g'+g''=g,n'+n''=n} \{ \mathcal{V}_{g',n'}^{[0]}, \mathcal{V}_{g'',n''}^{[0]}\}- \Delta \mathcal{V}_{g-1,n+2}^{[0]}
 \end{equation}
The boundary of the expression $q{\sf pr} \mathbb{I}i^\sharp \big(\mathcal{K}_*\mathcal{V}\big)_{g,n}$ consists of three types of terms:
\begin{itemize}
\item[(i)] For an edge $e\in E_G$, the corresponding parameter $t_e=0$ (with a minus sign).
\item[(ii)] For an edge $e\in E_G$, the corresponding parameter $t_e=1$ (with a plus sign).
\item[(iii)] At a vertex $v\in V_G$, there is $\partial\mathcal{V}_{g(v),n(v)}^{\sym,[0]} $.
\end{itemize}
We observe that using Equation~\eqref{eq:leading-mc}, the $(ii)$ and $(iii)$ type terms cancel each other. Moreover, the $(i)$ term is precisely the twisted sewing operations used in defining the DGLA structure of $\mathfrak{n}$. Hence depending on whether $e\in E_G$ is a separating or a non-separating edge, the $(i)$ term can be rewritten as
\[ - \frac{1}{2} \sum_{g'+g''=g,n'+n''=n} \{ q{\sf pr} \mathbb{I}i^\sharp \big(\mathcal{K}_*\mathcal{V}\big)_{g',n'}, q{\sf pr} \mathbb{I}i^\sharp \big(\mathcal{K}_*\mathcal{V}\big)_{g'',n''}\}- \Delta q{\sf pr} \mathbb{I}i^\sharp \big(\mathcal{K}_*\mathcal{V}\big)_{g-1,n+2}.\]
This finishes the proof.
\end{proof}

\section{CEI of the ground field}~\label{sec:cei}

In this section, we use Equation~\eqref{eq:fundamental} to prove Theorem A. Throughout the section, let $A$ be a cyclic $A_\infty$ algebra of dimension $d$ over a field $\mathbb{K}$ (of characteristic zero) that is proper, smooth and satisfies the Hodge-to-de-Rham degeneration.

\paragraph{{\bf A sketch of the definition of CEI.}} CEI defined in~\cite{CT}) takes a roundabout route due to the fact that the $2$-dimensional Topological Conformal Field Theory (TCFT) structure on the reduced Hochschild chain complex $C_*(A)$ has to have strictly positive number of inputs. This separation of inputs and outputs is the source of main difficulties in making Costello's original definition~\cite{Cos2} of categorical enumerative invariants explicit. Indeed, in~\cite{CCT} and~\cite{CT}, a Koszul type resolution is constructed for the Sen-Zwiebach Lie algebra (see Paragraph~\ref{para:DGLA}), giving a quasi-isomorphism of DGLA's:
\[ \mathfrak{g}  \ra   \widehat{\mathfrak{g}},\]
where the right hand side DGLA $ \widehat{\mathfrak{g}}$ consists of Riemann surfaces with strictly positive number of inputs. The precise definition of $ \widehat{\mathfrak{g}}$ is not so much relevant to this paper: we shall see that in the particular case when the algebra is the ground field there is no need to use the Koszul resolution.

Similar to the construction of the Sen-Zwiebach Lie algebra $\mathfrak{g}$, associated with the Hochschild chain complex $C_*(A)[d]$ and the Connes operator $B$, we may form a DGLA 
\[ \mathfrak{h}_A := \sym \big( C_*(A)[d]\big)[1] [[\hbar,\lambda]].\]
In~\cite{CCT}, we also introduced its Koszul resolution $\widehat{\mathfrak{h}}_A$, and a canonical quasi-isomorphism of DGLA's:
\[ \mathfrak{h}_A \ra  \widehat{\mathfrak{h}}_A.\]
Again we omit the definition of $\widehat{\mathfrak{h}}_A$ and refer the details to~\cite{CCT}. Note that both $\mathfrak{h}_A$ and $\widehat{\mathfrak{h}}_A$ are pro-nilpotent DGLA's in the $(\hbar,\lambda)$-adic topology. For such DGLA's, when considering its Maurer-Cartan elements, we shall always require they lie inside $(\hbar,\lambda)\cdot \mathfrak{h}_A$ or $(\hbar,\lambda)\cdot \widehat{\mathfrak{h}}_A$. The homotopy invariance of the associated Maurer-Cartan moduli space is proved in~\cite{Yek}.

In an ideal situation, the TCFT structure on $C_*(A)$ would give us a commutative diagram 
\[\begin{tikzpicture}[scale=1.5]
\node (A) at (0,1) {$\mathfrak{g}$};
\node (B) at (1,1) {$\widehat{\mathfrak{g}}$};
\node (C) at (0,0) {$\mathfrak{h}_A$};
\node (D) at (1,0) {$\widehat{\mathfrak{h}}_A$};
\path[->,font=\scriptsize,>=angle 90]
(A) edge node[above]{$\iota$} (B)
(B) edge node[right]{$\widehat{\rho}_A$} (D)
(C) edge node[above]{$\iota$} (D);
\path[->,dashed,font=\scriptsize,>=angle 90]
(A) edge node[right]{$\rho_A$} (C);
\end{tikzpicture}\]
with $\rho_A$ and $\widehat{\rho}_A$ the TCFT structure maps. In reality, the left vertical map $\rho_A$ is not available as a DGLA homomorphism, which prevents us to obtain the push-forward $(\rho_A)_*\mathcal{V}$ of the string vertex from $\mathfrak{g}$ to define categorical enumerative invariants. One might hope to construct $\rho_A$ as an $L_\infty$-morphism. The existence of such a morphism is clear from the diagram above since the bottom map $\iota$ admits an $L_\infty$-inverse. However, no explicit construction of such an $\rho_A$ is known yet.

In~\cite{CT}, we took the roundabout route:
\begin{itemize}
\item[(1.)] We find string vertices $\widehat{\mathcal{V}}$ in $\widehat{\mathfrak{g}}$.
\item[(2.)] Then push-forward $\widehat{\mathcal{V}}$ to obtain a Maurer-Cartan element $$\widehat{\beta}:=(\widehat{\rho}_A)_*\widehat{\mathcal{V}}$$ 
of the DGLA $\widehat{\mathfrak{h}}_A$.
\item[(3.)] The bottom map $\iota$ is a quasi-isomorphism of DGLA's when $A$ is smooth, proper, and satisfies the Hodge-to-de-Rham degeneration property. Hence, by homotopy invariance of the Maurer-Cartan moduli space of pro-nilpotent DGLA's, we obtain the desired Maurer-Cartan element $\beta\in \mathfrak{h}_A$ with its defining property that $\iota_*\beta$ is gauge-equivalent to $\widehat{\beta}$.
\end{itemize}

\paragraph{{\bf Categorical enumerative invariants of the ground field.}} When $A$ is the ground field $\mathbb{Q}$, the above roundabout route is not necessary: we can actually define a map
\[ \rho_A: \mathfrak{g}\ra \mathfrak{h}_A\]
to obtain the commutative diagram in the previous paragraph. Thus when computing the CEI, we may simply take the Maurer-Cartan element $\beta=(\rho_A)_*\mathcal{V}$ in $\mathfrak{h}_A$. Explicitly, the map $\rho_A$ is just the augmentation map: at degree $(g,n)$, the map $\rho_A: C_*(M_{g,n}^{\sf fr}) \ra \mathbb{Q}$ sends any point class in $C_0(M_{g,n}^{\sf fr})$ to $1$, and is zero otherwise.

For the splitting of the Hodge filtration, observe that the Connes operator $B$  acts on  by zero on $C_*(A)=\mathbb{Q}$. Thus, we may choose the splitting $s: C_*(A) = \mathbb{Q}\ra C_*(A)[[u]]=\mathbb{Q}[[u]]$ to be
\[ s(1)= 1.\]
Since $u$ is of degree $-2$, such a splitting is characterized by the homogeneity condition. By Proposition~\ref{prop:extension}, a splitting induces an extension of the TCFT to give a map
\[ \overline{\rho}_{A,s}: \overline{\mathfrak{g}} \ra \mathfrak{h}_A.\]
In summary, we obtain a commutative diagram of $L_\infty$ morphisms between DGLA's:
\begin{equation}~\label{diag:trivialization}
\begin{CD}
\mathfrak{g}@>i>> \overline{\mathfrak{g}} @>\mathcal{K} >> \overline{\mathfrak{g}}^{\sf triv}\\
@. @V\overline{\rho}_{A,s} VV                       @V\overline{\rho}_{A,s} VV\\
@. \mathfrak{h}_A   @> \mathcal{K}_{A,s}>> \mathfrak{h}_A^{\sf triv}
\end{CD}
\end{equation}
The construction of $\mathcal{K}$ is done in Section~\ref{sec:fundamental} while $\mathcal{K}_A$ is done in~\cite{CT}.

\begin{Theorem}~\label{thm:main2}
The CEI of $(A=\mathbb{Q},s)$ agrees with the Gromov-Witten invariants of a point, i.e. we have
\[ \langle \psi^{k_1}, \ldots, \psi^{k_n} \rangle_{g,n}^{\mathbb{Q},s} = \int_{\overline{M}_{g,n}} \psi_1^{k_1}\cdots \psi_n^{k_n}.\]
\end{Theorem}

\begin{proof}
The string vertex $\mathcal{V}_{g,n}$ is of the form
\[ \mathcal{V}_{g,n} = [\sum_{k_1,\ldots,k_n \geq 0} \mathcal{V}_{g,n}^{k_1,\ldots,k_n} u_1^{-k_1}\cdots u_n^{-k_n}],\]
with $ \mathcal{V}_{g,n}^{k_1,\ldots,k_n} \in C_{6g-6+2n-2k_1-\cdots -2k_n}(M_{g,n}^{\sf fr}) $. The bracket is to indicate the quotient map by the $S_n$ action. In characteristic zero, we may assume that the chains $ \mathcal{V}_{g,n}^{k_1,\ldots,k_n}$ are invariant under the action of $S_n$. Then we have
\[ \overline{\rho}_{A,s} i (\mathcal{V}_{g,n}) = \sum_{k_1+\cdots+k_n=3g-3+n} \rho_A( \mathcal{V}_{g,n}^{k_1,\ldots,k_n})\frac{n!}{{\sf Aut}(k_1,\ldots,k_n)} q_{k_1}\cdots q_{k_n}.\]
In the above formula, the right hand side takes value in $\mathfrak{h}_A= \sym (\mathbb{Q}[u^{-1}])$, and the variable $q_{k_j}$ corresponds the basis $u^{-k_j}$. The combinatorial number ${\sf Aut}(k_1,\ldots,k_n)$ is the size of the set of $n$-permutations that fixes the sequence $(k_1,\ldots,k_n)$. The extension formula  in Proposition~\ref{prop:extension}  of the splitting $s$ yields the zero map for the action of $D_e$ on the edges of a stable graph. This implies that the bottom $L_\infty$ algebra map $\mathcal{K}_{A,s}$ is the identity map, which shows that
\[\mathcal{K}_{A,s} \overline{\rho}_{A,s} i (\mathcal{V}_{g,n}) = \sum_{k_1+\cdots+k_n=3g-3+n} \rho_A( \mathcal{V}_{g,n}^{k_1,\ldots,k_n})\frac{n!}{{\sf Aut}(k_1,\ldots,k_n)} q_{k_1}\cdots q_{k_n}.\]
By definition of CEI, we have
\begin{align*} 
& \langle \psi^{l_1}, \ldots, \psi^{l_n} \rangle_{g,n}^{\mathbb{Q},s} \\
= & \frac{\partial}{\partial q_{l_1}}\cdots\frac{\partial}{\partial q_{l_n}}\big( \sum_{k_1+\cdots+k_n=3g-3+n} \rho_A( \mathcal{V}_{g,n}^{k_1,\ldots,k_n})\frac{n!}{{\sf Aut}(k_1,\ldots,k_n)} q_{k_1}\cdots q_{k_n}\big)\\
= & n! \rho_A( \mathcal{V}_{g,n}^{l_1,\ldots,l_n})
\end{align*}
On the other hand, by Theorem~\ref{thm:main} we have
\[ [\overline{M}_{g,n}/S_n]= (\mathcal{K}_*\mathcal{V})_{g,n}= \sum_{G\in \Gamma((g,n))} \frac{1}{|{\sf Aut}(G)|} \prod_{e\in E_G} D_e \otimes \prod_{v\in V_G} \mathcal{V}^{\sym}_{g(v),|{\sf Leg}(v)|}.\]
This implies that
\begin{align*}
& \int_{\overline{M}_{g,n}} \psi_1^{l_1}\cdots \psi_n^{l_n}\\
= & n!\cdot  u_1^{l_1}\cdots u_k^{l_n} \big( \sum_{G\in \Gamma((g,n))} \frac{1}{|{\sf Aut}(G)|} \prod_{e\in E_G} D_e \otimes \prod_{v\in V_G} \mathcal{V}^{\sym}_{g(v),|{\sf Leg}(v)|} \big)\\
= & n! \cdot u_1^{l_1}\cdots u_k^{l_n} \big( \sum_{k_1,\ldots,k_n \geq 0} \mathcal{V}_{g,n}^{k_1,\ldots,k_n} u_1^{-k_1}\cdots u_n^{-k_n}\big)\\
=& n! \rho_A( \mathcal{V}_{g,n}^{l_1,\ldots,l_n})
\end{align*}
This proves the theorem.
\end{proof}

\end{document}